\documentclass{amsart}

\usepackage{amsmath}
\usepackage{amssymb}
\usepackage{latexsym,amscd}
\input V3,97

\title[Narrow  towers]
          {  On the Narrow $2$-Class Field Tower of  \\ Some Real Quadratic Number Fields \\ Part II: Lengths }

\author{E.\,Benjamin,  C.\,Snyder}

\keywords{real quadratic field, Hilbert class field, 2-class group, narrow class group, narrow 2-class field tower}
\subjclass[2010]{11R29, 20D15}

\begin{document}
\newcommand{\rsp}{\raisebox{0em}[2.7ex][1.3ex]{\rule{0em}{2ex} }}
\newcommand{\ul}{\underline}
\newcommand{\1}{\boldsymbol{1}}
\newcommand{\I}{\boldsymbol{i}}
\newcommand{\J}{\boldsymbol{j}}
\newcommand{\K}{\boldsymbol{k}}
\newcommand{\C}{{\mathbb C}}
\newcommand{\h}{\mathbb{H}}
\newcommand{\Q}{{\mathbb Q}}
\newcommand{\R}{{\mathbb R}}
\newcommand{\B}{{\mathbb B}}
\newcommand{\F}{{\mathbb F}}
\newcommand{\N}{{\mathbb N}}
\newcommand{\Z}{{\mathbb Z}}
\newcommand{\aQ}{{\overline{\Q}}}
\newcommand{\cA}{\mathcal A}
\newcommand{\cB}{\mathcal B}
\newcommand{\cF}{\mathcal F}
\newcommand{\cN}{\mathcal N}
\newcommand{\cK}{\mathcal K}
\newcommand{\OO}{\mathcal O}
\newcommand{\cP}{\mathcal P}
\newcommand{\Rl}{{\rm Re}}
\newcommand{\cM}{\mathcal M}
\newcommand{\Am}{\operatorname{Am}}
\newcommand{\eq}{\stackrel{2}{=}}
\newcommand{\fZ}{\mathfrak 2}
\newcommand{\fQ}{\mathfrak Q}
\newcommand{\fa}{\mathfrak a}
\newcommand{\fb}{\mathfrak b}
\newcommand{\fc}{\mathfrak c}
\newcommand{\fC}{\mathfrak C}
\newcommand{\fE}{\mathfrak E}
\newcommand{\fp}{\mathfrak p}
\newcommand{\fP}{\mathfrak P}
\newcommand{\cO}{\mathcal O}
\newcommand{\fA}{\mathfrak A}
\newcommand{\fB}{\mathfrak B}
\newcommand{\fG}{\mathfrak G}
\newcommand{\fr}{\mathfrak r}
\newcommand{\frb}{\widetilde{\fr}}
\newcommand{\fR}{\mathfrak R}
\newcommand{\fq}{\mathfrak q}
\newcommand{\gac}{\widehat}
\newcommand{\Ft}{\widetilde{F}}
\newcommand{\disc}{\operatorname{disc}}
\newcommand{\ch}{\operatorname{char}}
\newcommand{\id}{\operatorname{id}}
\newcommand{\rank}{\operatorname{rank}}
\newcommand{\lcm}{{\operatorname{lcm}}}
\newcommand{\gen}{{\operatorname{gen}}}
\newcommand{\fin}{{\operatorname{fin}}}
\newcommand{\cen}{{\operatorname{cen}}}
\newcommand{\sfk}{{\operatorname{sfk}}}
\newcommand{\Tr}{\operatorname{Tr}}
\newcommand{\Cl}{\operatorname{Cl}}
\newcommand{\Ram}{\operatorname{Ram}}
\newcommand{\Gal}{\operatorname{Gal}}
\newcommand{\Aut}{\operatorname{Aut}}
\newcommand{\AGL}{\operatorname{AGL}}
\newcommand{\Ind}{\operatorname{Ind}}
\newcommand{\GL}{\operatorname{GL}}
\newcommand{\im}{\operatorname{im}\,}
\newcommand{\eps}{\varepsilon}
\newcommand{\la}{\langle}
\newcommand{\ra}{\rangle}
\newcommand{\lra}{\rightarrow}
\newcommand{\too}{\mapsto}
\newcommand{\ts}{\textstyle}
\newcommand{\ab}{\operatorname{ab}}
\newcommand{\ff}{\mathfrak f}
\newcommand{\fm}{\mathfrak m}
\newcommand{\fn}{\mathfrak n}
\def\la{\langle}
\def\ra{\rangle}
\def\lra{\rightarrow}
\def\lms{\mapsto}
\def\vp{\varphi}
\def\vt{\vartheta}
\def\tr{\mbox{{\rm Tr}}\,}
\def\hom{\mbox{{\rm Hom}}\,}
\def\ind{\mbox{{\rm Ind}}}
\newtheorem{df}{Definition}
\newtheorem{thm}{Theorem}
\newtheorem{lem}{Lemma}
\newtheorem{prop}{Proposition}
\newtheorem{cor}{Corollary}

\begin{abstract} We determine precisely when the length of the  narrow 2-class field tower is $2$ for most of those real quadratic number fields whose  discriminant is not a sum of two squares and for which their 2-class groups are elementary of order $4$.
\end{abstract}

\maketitle

\section{{\bf Introduction}}

This is a continuation of an investigation into the narrow $2$-class field tower of those real quadratic fields $k$ whose  discriminants are not a sum of two squares and for which their 2-class groups are elementary of order $4$;  see \cite{BS3} for relevant results and notation used here. For instance, let $k^i$ and $k_+^i$ denote the $i$-th ordinary, resp., narrow $2$-class field of $k$; moreover, let $G=\Gal(k^2/k)$ and $G^+=\Gal(k_+^2/k)$. We are now interested in determining the length of the narrow $2$-class field tower of these fields $k$,  and  actually we have made some progress already.  Namely, for all those $k$ for which $G^+/G^+_3$ are of order $32$, but not of type $32.033$ (see \cite{HS} for group designations), the length of the tower is $2$, by Corollary~1 and Tables 1 and $2$ in \cite{BS3}. In particular,  $G^+/G^+_2$ is elementary and $G^+_2$ has rank $2$, which implies that $G^+$ is metabelian (refer to the Main Theorem in \cite{BS2}).

On the other hand, we now investigate the length of the tower when $G^+_2$ has rank at least $3$. We establish the following main result:

\bigskip
\noindent {\bf Main Theorem.} {\em Let $k$ be any real quadratic field whose  discriminant is not a sum of two squares and for which its $2$-class group is elementary of order $4$. If  either $G^+/G_3^+\simeq 32.033$ or of order $64$, with the exception of  group $64.150$, then  the narrow $2$-class field tower of $k$ has length $\geq 3$.}

\bigskip

For the case $G^+/G^+_3\simeq 64.150$, all our examples have narrow $2$-class field tower of length $2$, however, the general case escapes us. Perhaps there is an example of length at least $3$, but we have not been able to exhibit any. We have  been able to show that if the $8$-rank of $\Cl_2(k_+^1)$ is $0$, i.e., if $\Cl_2(k_+^1)\subseteq (4,4,4)$, then the narrow $2$-class field tower has length $2$.

We'll prove the theorem  by constructing an  unramified quadratic extension $L$ of $k_+^1$ and showing that $h_2(L)/h_2(k_+^1)\geq 1.$
More precisely, in all these cases (including $64.150$) it turns out that there is a quadratic subfield $F$ of $k_+^1$ for which its discriminant has a $C_4$-splitting.  From the $C_4$-splitting we obtain a finitely unramified (i.e.,\,unramified outside $\infty$) cyclic quartic extension $F_c$ of $F$. If we let $L=k_+^1F_c$, then we have an unramified quadratic extension of $k_+^1$ which  is a $V_4$-extension of $k^1$.  We then use Kuroda's class number formula to get a handle on the quotient $h_2(L)/h_2(k_+^1)$, by considering arithmetic properties of the intermediate fields in the extension $L/k^1$.   We will consider four families of cases:

\medskip
(A) the cases for which $G$ is either dihedral, semi-dihedral, or generalized quaternion,

(B) the cases for which $G$ is quaternion,

(C) the cases for which $G$ is abelian, except for $G^+/G^+_3\simeq 64.150,$

(D) the cases for which $G^+/G^+_3\simeq 64.150.$
\medskip

Each family will require somewhat different methods.

\section{{\bf (A) the cases for which $G$ is $D,S,$ or $Q_g$}}

We'll start by considering the cases of $(a_i), (b_i)$ for $i=1,\dots 8$ as given in Appendix~I below. (We have included the case $(b_1)$ at no real extra cost  even though we already know the length of the narrow $2$-class field tower of $k$ is $2$.) These are the cases where $G=\Gal(k^2/k)$ is dihedral, semi-dihedral, or generalized quaternion. Appendix~I consolidates some useful information from the appendices in \cite{BS3}.

Since $k$ satisfies any  of these cases, we see by Appendix~I (or II)  that there is a quadratic subfield $F$ of $k^1$ whose discriminant has a $C_4$-splitting. In fact, let
$$F=\Q(\sqrt{d_1D}\,)\;\text{such that} \;\;D=\left\{
                                            \begin{array}{ll}
                                              d_2, & \;\hbox{if}\;\; k\in A(d_2), \\
                                              d_3d_4, & \;\hbox{if}\;\; k\in A(d_3d_4),
                                            \end{array}
                                          \right. $$
where $A(D)$ is the set of cases such that
\begin{align*}
 A(d_2)& =\{(a_i), (b_i): i=1,2,5,6,7\},\\
A(d_3d_4)& =\{ (a_i),(b_i): i=3,4,8\}.
\end{align*}
Thus the discriminant of $F$ satisfies the $C_4$-splitting $d_F=d_1\cdot D$,
 and hence  $\Cl_2^+(F)\supseteq C_4$. Observe that  $\Cl_2(F)\supseteq C_4$, if $h_2(d_1D)\geq 4$; however, when $h_2(d_1D)=2$, then $\Cl_2^+(F)\simeq C_4$;  moreover this can only occur when $D=d_2$. In both instances let $F_c$ be the intermediate field such that $F\subseteq F_c\subseteq F_+^1$ with $F_c/F$ a cyclic quartic extension unramified at least at all the finite primes. Notice that $F_c$ is a quadratic extension of  $F_\gen=\Q(\sqrt{d_1},\sqrt{D}\,)$. Now, we have $F_c=F(\sqrt{\alpha}\,)$, where for a particular rational solution $(a,b,c)\in\Q^3$ of $x^2=y^2d_1+z^2D$, $\alpha=a+c\sqrt{D}$, cf.\,\cite{Lem}. We let $\alpha'=a-c\sqrt{D}$; then $N\alpha=\alpha\alpha'=a^2-c^2D=b^2d_1$. Hence $\sqrt{\alpha'}=\pm b\sqrt{d_1}/\sqrt{\alpha}\in F_c$; but $\Q(\sqrt{\alpha}\,)$ is not normal over $\Q$. Observe that for $h_2(F)\geq 4$, $\alpha$ and $\alpha'$ are both positive, since $F_c$ is totally real. On the other hand, when $h_2(F)=2$, then both $\alpha$ and $\alpha'$ are negative. For use below we let
$$\gamma=\left\{
           \begin{array}{ll}
             \alpha, & \hbox{if}\;\; \alpha<0, \\
             d_3\alpha, & \hbox{if}\;\; \alpha>0.
           \end{array}
         \right. $$

Next, let $L=k_+^1(\sqrt{\alpha}\,)$. This is a quadratic extension of $k_+^1$, since $\alpha\in k_+^1$ but $\sqrt{\alpha}$ is not, because $k_+^1/\Q$ is abelian, but $\Q(\sqrt{\alpha}\,)/\Q$ is not normal. Since $F_c/F$ is a cyclic, quartic extension unramified outside $\infty$, the same is true of $kF_c/kF$ (notice that it remains quartic, since $F_c\cap kF=F$). Also, since $kF/k$ is unramified everywhere, $kF_c/k$ is finitely unramified of degree $8$. Moreover, $L=kF_c(\sqrt{d_3})$ and thus $L/k$ is finitely unramified which in turn implies that
$k_+^1\subseteq L\subseteq k_+^2.$  We see that $k^1=kF_\gen$ and that $L=k^1(\sqrt{\alpha},\sqrt{d_3}\,)$. Thus $L/k^1$ is a $V_4$-extension with intermediate fields $k_+^1,\, K=k^1(\sqrt{d_3\gamma}\,)$, and $N=k^1(\sqrt{\gamma}\,)$. We have

$$ L=k^1(\sqrt{d_3},\sqrt{\gamma}\,),$$
$$ k_+^1=k^1(\sqrt{d_3}\,),\;\;\;K=k^1(\sqrt{d_3\gamma}\,),\;\;\; N=k^1(\sqrt{\gamma}\,),$$
$$ k^1=\Q(\sqrt{d_1},\sqrt{d_2},\sqrt{d_3d_4}\,).$$

We observe that $k_+^1,\, N,\,L $ are totally complex, whereas $K$ is totally real. In fact $k_+^1,\, N,\, L$ are CM-fields with maximal real subfields $k^1,\,k^1,\,K$, respectively. Our goal is to show that
$h_2(L)/h_2(k_+^1)\geq 1$, except for the case $(b_1)$ where we already know $h_2(L)/h_2(k_+^1)=1/2$. We will be applying Kuroda's class number formula several times to fields appearing in the following diagram:

\bigskip
Diagram 1

\bigskip
\begin{diagram}[height=2em,width=2em]
& &  & & & &  L & &  \\
 & & & & & \ldLine & \dLine &\rdLine  &  \\
 & & & & N & & K & & k_+^1\\
 & &  &\ldLine &\dLine &\rdLine & \dLine &\ldLine  &  \\
 &  & M & & M'  & & k^1 & &  \\
 & \ldLine  & \dLine &\rdLine & \dLine & \ldLine & & &  \\
\Q(\sqrt{D_1\gamma}\,) & & \Q(\sqrt{\gamma}\,) & & \mathcal{F} & & & & \\
& \rdLine &  \dLine  &\ldLine &\dLine &\rdLine &  &  &  \\
 & & \Q(\sqrt{D}) & & \Q(\sqrt{d_2d_3d_4}) & & \Q(\sqrt{D_1})  &  &\\
 & &  &\rdLine & \dLine & \ldLine & & &  \\
 &&&& \Q &&&& \\
 \end{diagram}

\noindent where

$$ D_1=d_2d_3d_4/D $$
$$\mathcal{F}=\Q(\sqrt{d_2},\sqrt{d_3d_4}\,),$$
$$ M=\mathcal{F}(\sqrt{\gamma}\,),\; M'=\mathcal{F}(\sqrt{\gamma'}\,).$$

Roughly speaking, we will use Kuroda's class number formula several times to reduce computing lower bounds on $h_2(L)/h_2(k_+^1)$ to computing lower bounds on the $2$-class numbers of the two quartic fields $\Q(\sqrt{\gamma}\,)$ and $\Q(\sqrt{D_1\gamma}\,)$.  More precisely, we will show the following results:

\begin{prop}\label{PK1} Let $\ell=\left\{
                                    \begin{array}{ll}
                                      1/4, & \hbox{if}\; G\simeq Q_g, \\
                                      1/8, & \hbox{if}\; G\simeq D, S.
                                    \end{array}
                                  \right.$
Then by Kuroda's class number formula applied to $L/k^1$,
$$\frac{h_2(L)}{h_2(k_+^1)}=\ell\,\frac{h_2(N)}{h_2(k^1)}.$$
\end{prop}

\begin{prop}\label{PK2} By Kuroda's class number formula applied to $\mathcal{F}/\Q$,
$$h_2(\mathcal{F})=1.$$
\end{prop}

\begin{prop}\label{PK3} By Kuroda's class number formula applied to $N/\mathcal{F}$,
$$\frac{h_2(N)}{h_2(k^1)}=\frac1{8\,\ell}\,h_2(M)^2.$$
\end{prop}

Therefore
$$\frac{h_2(L)}{h_2(k_+^1)}=\frac1{8}\,h_2(M)^2.$$
Finally,
\begin{prop}\label{PK4} By Kuroda's class number formula applied to $M/\Q(\sqrt{D})$,
$$ h_2(M)=\frac1{2} h_2\big(\Q(\sqrt{\gamma})\big)h_2\big(\Q(\sqrt{D_1\gamma})\big).$$
\end{prop}

Therefore
$$\frac{h_2(L)}{h_2(k_+^1)}=\frac1{32}h_2\big(\Q(\sqrt{\gamma})\big)^2h_2\big(\Q(\sqrt{D_1\gamma})\big)^2.$$

\bigskip
{\em We now start with the proof of Proposition~\ref{PK1}.}

\medskip
Applying Kuroda's class number formula (cf.\,\cite{Lem2}) to $L/k^1$  yields
$$h_2(L)/h_2(k_+^1)=\frac1{2^{1+v}}\,q(L/k^1) h_2(K)\,h_2(N)/h_2(k^1)^2$$
(since $d=8$, $\kappa=7$, and $v=0$ or $1$).  But by \cite{Lem3}, Proposition~3, since $L/k^1$ is a $V_4$-extension of CM-fields,
$$q(L/k^1)=2^{1+v}\,\frac{Q(L)}{Q(k_+^1)Q(N)}\,\frac{w_L}{w_{k_+^1}w_N},$$
where $w$ is the order of the group $W$ of the roots of unity in the relevant field and $Q(L)=(E_L:W_L E_K)$, etc.

\begin{lem}\label{L3.1} For $L, k_+^1,N$ as above, $w_L/(w_{k_+^1}w_N)=1/2$.
\end{lem}
\begin{proof} Since all ramification at finite primes in $L/\Q$ occurs in $k/\Q$, the ramification indices are at most $2$. Hence the only possible roots of unity in $L$ are those contained in $\Q(\sqrt{d_3},\sqrt{d_4}\,)$. Clearly $W_{k_+^1}\subseteq W_L$; but the opposite inclusion is also clear since $\Q(\sqrt{d_3},\sqrt{d_4}\,)\subseteq k_+^1.$ Thus $w_L=w_{k_+^1}.$ On the other hand, we observe that the only quadratic subfields contained in $N$ are all real, since otherwise $k_+^1\subseteq N$. Therefore, $w_N=2$. This establishes the lemma.
\end{proof}

We also see that $h_2(K)=h_2(k^1)/2$ since $\Cl_2(K)\simeq \Gal(K^1/K)=\Gal(k^2/K).$ In light of this fact and the lemma, we have
$$\frac{h_2(L)}{h_2(k_+^1)}=\frac1{4}\,\frac{Q(L)}{Q(k_+^1)Q(N)}\,\frac{h_2(N)}{h_2(k^1)}.$$

We now determine the fraction $Q(L)/Q(k_+^1)Q(N)$.

\begin{lem}\label{L3.2} The indices $Q(k_+^1)=2=Q(L)$.
\end{lem}
\begin{proof} First, suppose $w_{k_+^1}\equiv 2\bmod 4$. Hence $-4$ does not divide $d_k$. Notice that $k_+^1=k^1(\sqrt{d_3}\,)$ is not essentially ramified over $k^1$ (since $p_3$ ramifies in $k$).
Now $p_3\cO_{F_{34}}=(\pi_3^2)$ for a prime element $\pi_3$ in $F_{34}=\Q(\sqrt{d_3d_4}\,)$, since $p_3$ is ramified and $h_2(d_3d_4)=1$. Hence $p_3\cO_{k^1}=(\pi_3\cO_{k^1})^2$ and thus by Theorem~1(i) in \cite{Lem3}, $Q(k_+^1)=2$.
Similarly, $p_3\cO_K=(\pi_3\cO_K)^2$ and so $Q(L)=2.$

Now suppose $w_{k_+^1}\equiv 4\bmod 8$. Then $2\cO_{k^1}=\fq^2\cO_{k^1}$ for the  ideal $\fq$ of $k$ dividing $2$. But  by referring to Appendix II, notice that $\fq$ is either principal (in $(b_2), (b_3)$ since $\delta=2$ or $2p_1p_2p_3$) or capitulates in $k_3$ (in the rest of the cases). Therefore, $\fq\cO_{k^1}$ is a principal ideal and by Theorem~1(ii) in \cite{Lem3}, $Q(k_+^1)=2$. Similarly, $Q(L)=2$.
\end{proof}

\begin{lem}\label{L3.3} The index $Q(N)$ satisfies
$$Q(N)=\left\{
         \begin{array}{ll}
           1, & \hbox{if $G$ is generalized quaternion,} \\
           2, & \hbox{if $G$ is dihedral or semi-dihedral.}
         \end{array}
       \right.$$
\end{lem}
\begin{proof} First, assume $G$ is generalized quaternion. We claim that $N=k^1(\sqrt{d_3\alpha}\,)$ with $\alpha>0$. For, if not, then $N=k^1(\sqrt{\alpha}\,)$ with $\alpha<0$. In this case, $h_2(d_1D)=2$. Thus $D=d_2$ and so by Appendix II (cases $(a_1),(a_2),(b_1),(b_2)$)\; $|G|=4h_2(d_1d_2)=8$, contrary to our assumption that $G$ is generalized quaternion.

Hence assume $N=k^1(\sqrt{d_3\alpha}\,)$, where $\alpha\in F_D=\Q(\sqrt{D}\,)$ (as above) with $\alpha>0$. Since $N/k^1$ is unramified outside $\infty$, we see that $N/k^1$ is not essentially ramified and thus the ideal  $d_3\alpha\cO_{k^1}$ is an ideal  square  in $\cO_{k^1}.$ We have already seen from above that $p_3\cO_{k^1}=(\pi_3\cO_{k^1})^2$ (clearly then $d_3$ generates the square of a principal ideal, too). So now consider $\alpha\cO_{k^1}$. As we have seen, $N\alpha=\alpha\alpha'=a^2-c^2D=b^2d_1=\tilde{b}^2p_1$ (in case $d_1=8$). Thus $\alpha\cO_{F_D}=\tilde{b}\tilde{\fp}_1$, with $\tilde{\fp}_1$ a prime ideal. In fact, $p_1\cO_{F_D}=\tilde{\fp}_1\tilde{\fp}_1'$ where $\tilde{b}\tilde{\fp}_1'=(\alpha')$. Now extend this ideal to $k(\sqrt{D}\,)$. Since $p_1$ is ramified in $k$ but splits in $F_D$, we see that $p_1\cO_{k(\sqrt{D}\,)}=\fP_1^2\fP_1'^2$ for prime ideals $\fP_1,\fP_1'$. Thus $\alpha\cO_{k(\sqrt{D}\,)}=\tilde{b}\fP_1^2$, say. This implies that $\alpha\cO_{k^1}=\tilde{b}(\fP_1\cO_{k^1})^2$. But we now see that $\tilde{b}\cO_{k^1}$ is an ideal square in $\cO_{k^1}$. This implies  that $\tilde{b}\cO_k=\fb^2$ (because  any  rational prime that  divides $\tilde{b}$ exactly to an odd power must be ramified in $k^1/\Q$, hence already in $k$.) But since $\fb^2\sim 1$, we see that $\fb$ becomes principal in the $2$-class field $k^1$ of $k$. We thus have $d_3\alpha\cO_{k^1}=(\beta\fP_1)^2$ for a particular $\beta\in \cO_{k^1}.$

It suffices to show that $\fP_1\cO_{k^1}$ is not a principal ideal. Toward this end, notice that $\fP_1$ is inert in $k^1$; otherwise the prime ideal $\fp_1$ in $k$ above $p_1$ would split completely, which cannot happen since $\fp_1$ is not principal (cf.\,Appendix II). Let $\cP_1=\fP_1\cO_{k^1}$. Similarly we see that  $\fP_1'\cO_{k^1}=\cP_1'$ is a prime ideal.  Now, by Theorem~1 in \cite{Lem3}, to show that $Q(N)=1$, it suffices to show that the ideal $\cP_1$ is not principal.

By way of contradiction, assume that $\cP_1$ is principal. We will then show that the capitulation kernel $\ker j_{k^1/k_D}$ with $k_D=k(\sqrt{D}\,)$ has order $4$. This will lead to the contradiction since for $G$ generalized quaternion, for each of the three quadratic fields $k_1,k_2,k_3$, their capitulation kernels $\ker j_{k_i/k}$ have order $2$; see section~4 in \cite{BS3}, for example.

First consider $D=d_2$, in which case $k_D=k_2$ (these are the cases\\ $(a_1),(a_2),(b_1),(b_2)$ in Appendix II). Then by Appendix II, $\Cl_2(k)=\la [\fp_1],[\fp_2]\ra$;  $\Cl_2(k_1)\simeq\Cl_2(k_2)\simeq (2,2)$ (since $\Gal(k^2/k_i)$ is quaternion ($i=1,2$)\;) and $\Cl_2(k_3)$ is cyclic; moreover $\kappa_1=\ker j_{k_1/k}=\la [\fp_1]\ra$, $\kappa_2=\la [\fp_2]\ra$, and $\kappa_3=\la [\fp_1\fp_2]\ra.$  Since $[\fp_1]\not\in \kappa_2$, we see that $\fp_1\cO_{k_2}=\fP_1\fP_1'\not\sim 1.$ Hence $\Cl_2(k_2)=\{ 1, [\fP_1],[\fP_1'],[\fP_1\fP_1']\}.$ But since $[\fp_1]\in \kappa_1$, $\cP_1\cP_1'=\fp_1\cO_{k^1}$ is principal. Thus $\fP_1\fP_1'\cO_{k^1}=\cP_1\cP_1'\sim 1$ in $k^1$ and so $[\fP_1\fP_1']\in \ker j_{k^1/k_2}$. But by assumption $[\fP_1]\in \ker j_{k^1/k_2}$. Therefore, $|\ker j_{k^1/k_2}|=4$, a contradiction.

A similar argument works for $D=d_3d_4$.

Therefore, $Q(N)=1$ for $G$ generalized quaternion.

\bigskip
Next,  suppose $G$ is dihedral. Then $N=k^1(\sqrt{\alpha}\,)$ or $k^1(\sqrt{d_3\alpha}\,)$. In either case, we are reduced to showing that $\fP_1$ (as above) a prime ideal in $k_D$ above $p_1$ becomes principal when extended to $k^1.$ We observe from Appendix II that $\Cl_2(k_D)\simeq (2,2)$. Thus since $G$ is dihedral, $\Gal(k^2/k_D)$ is either dihedral or $V_4$ (if $|G|=8$). Moreover, since  $\Gal(k^2/k^1)$ is the maximal cyclic proper subgroup of $\Gal(k^2/k_D)$, we see that $\ker j_{k^1/k_D}=\Cl_2(k_D).$ In particular $\fP_1\cO_{k^1}\sim 1$, as desired. Hence $Q(N)=2$.

\bigskip
Finally, assume $G$ is semi-dihedral. By Appendix II, we are looking at cases $(a_8), (b_8)$. In these cases, $D=d_3d_4$ and $k_D=k(\sqrt{D}\,)=k_3=k(\sqrt{d_1d_2}\,)$. Once again, we want to show that $\fP_1$ (as above) becomes principal in $k^1.$
By Appendix II again, we can see that $\Cl_2(k_2)\simeq \Gal(k^2/k_2)$ is cyclic, whereas $\Cl_2(k_i)\simeq (2,2)$ for $i=1, 3$. Since $G$ is semi-dihedral, one of the two maximal subgroups $\Gal(k^2/k_i)$ is quaternion, the other dihedral, ($i=1,3$). Since $[\fp_1]$ capitulates in both $k_1$ and $k_2$, but $[\fp_3]$ or $[\fq]$ capitulates in $k_3$, a group theoretic result implies that $\Gal(k^2/k_3)$ is dihedral, see e.g.\:p.\,131 in \cite{Lem}. (Also see Theorem~3 (3) in \cite{BS} and note  the remark/corrections in Section~4 of \cite{BS3}.)  Moreover, $\Gal(k^2/k^1)$ is the maximal cyclic subgroup of $\Gal(k^2/k_3)$. Thus $\ker j_{k^1/k_3}=\Cl_2(k_3)$ and this in turn implies in particular that $\fP_1$ becomes principal in $k^1$. Therefore, $Q(N)=2.$
\end{proof}

 Then
$$\frac{h_2(L)}{h_2(k_+^1)}=\frac{h_2(N)}{h_2(k^1)}\cdot \left\{
                                                           \begin{array}{ll}
                                                             1/4, & \hbox{if $G$ is generalize quaternion,} \\
                                                             1/8, & \hbox{if $G$ is dihedral or semi-dihedral.}
                                                           \end{array}
                                                         \right.
$$
This establishes Proposition~\ref{PK1}.

\bigskip
{\em Now we prove Proposition~\ref{PK2}}.

\medskip
We apply Kuroda's class number formula to the extension $\mathcal{F}/\Q$, cf.\;the diagram above. Then
$$h_2(\mathcal{F})=\frac1{4}q(\mathcal{F}/\Q)h_2(d_2)h_2(d_2d_3d_4)h_2(d_3d_4)=\frac1{4}q(\mathcal{F}/\Q)h_2(d_2d_3d_4).$$
But $q(\mathcal{F}/\Q)\leq 2$ and $h_2(d_2d_3d_4)=2$, since $\Cl_2(d_2d_3d_4)$ is cyclic and the discriminant of $\Q(\sqrt{d_2d_3d_4}\,)$ has no $C_4$-splittings. Thus $h_2(\mathcal{F})=1$, which proves Proposition~\ref{PK2}.

\bigskip
{\em Next, we prove Proposition~\ref{PK3}.}

\medskip
Applying Kuroda's class number formula to $N/\mathcal{F}$ yields
$$\frac{h_2(N)}{h_2(k^1)}=\frac1{2^{1+v}}q(N/\mathcal{F})h_2(M)h_2(M')/h_2(\mathcal{F})^2
=\frac1{2^{1+v}}q(N/\mathcal{F})h_2(M)^2.$$
Since $N,M,M'$ are CM-fields,\; $q(N/\mathcal{F})=2^vQ(N)/Q(M)^2$. But $Q(M)=1$, since $M/\mathcal{F}$ is essentially ramified. This implies that
$$\frac{h_2(N)}{h_2(k^1)}=\frac1{2}Q(N)h_2(M)^2=\left\{
                                                  \begin{array}{ll}
                                                   \frac1{2} h_2(M)^2, & \hbox{if}\;\;G\simeq Q_g, \\
                                                   \phantom{\frac1{2}} h_2(M)^2, & \hbox{if}\;\;G\simeq D, S.
                                                  \end{array}
                                                \right. $$
This proves Proposition~\ref{PK3}.

\bigskip
As noted above, Propositions~\ref{PK1} and \ref{PK3} imply that
$$\frac{h_2(L)}{h_2(k_+^1)}=\frac1{8}h_2(M)^2.$$
Also observe that
$$\frac{h_2(L)}{h_2(k_+^1)}\left\{
                          \begin{array}{ll}
                            =1/2, & \hbox{if}\;\;h_2(M)=2, \\
                            \geq 2, & \hbox{if}\;\; h_2(M)\geq 4.
                          \end{array}
                        \right.$$
Therefore to show that the narrow $2$-class field tower of $k$ has length $\geq 3$  for all cases except $(b_1)$, we need only show that $h_2(M)\geq 4$.

\bigskip
{\em We now prove Proposition~\ref{PK4}.}

\medskip
Applying Kuroda's class number formula to $M/\Q(\sqrt{D}\,)$, which is a $V_4$-extension of CM-fields, we get
$$h_2(M)=\frac1{2}\frac{Q(M)}{Q(\Q(\sqrt{\gamma}\,))Q(\Q(\sqrt{D_1\gamma}\,))}h_2(\mathcal{F})h_2(\Q(\sqrt{\gamma}\,))
h_2(\Q(\sqrt{D_1\gamma}\,))/h_2(D)^2.$$
But $h_2(\mathcal{F})=1$, $h_2(D)=1$, and $Q(M)=1$; moreover $Q(\Q(\sqrt{\gamma}\,))=Q(\Q(\sqrt{D_1\gamma}\,))=1$ since both $\Q(\sqrt{\gamma}\,)$ and $\Q(\sqrt{D_1\gamma}\,)$ are essentially ramified over $\Q(\sqrt{D}\,).$ Therefore
$$h_2(M)=\frac1{2}h_2(\Q(\sqrt{\gamma}\,))h_2(\Q(\sqrt{D_1\gamma}\,)),$$
proving Proposition~\ref{PK4}.

\bigskip
Now the real work comes in getting a handle on $h_2(\Q(\sqrt{\gamma}\,))$ and $h_2(\Q(\sqrt{D_1\gamma}\,))$. We'll consider the two cases $D=d_2$ and $D=d_3d_4$ separately (except for the next lemma).
\begin{lem}\label{L4} For $D\in\{d_2,d_3d_4\}$,
$$h_2(\Q(\sqrt{\alpha}\,))=1.$$
\end{lem}
\begin{proof} Consider the diagram:
\bigskip

Diagram 2

\bigskip
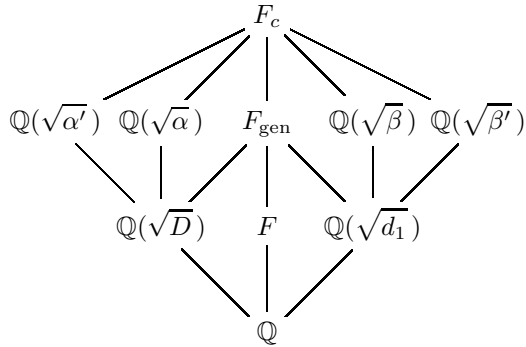
\begin{diagram}[height=2em,width=2em]
&   &   &    &  F_c   &   &   &   & \\
&   & & \ldLine(4,2)  \ldLine   &  \dLine             & \rdLine   \rdLine(4,2)  &   & & \\
\Q(\sqrt{\alpha'}\,) &   & \Q(\sqrt{\alpha}\,) & & F_\gen & & \Q(\sqrt{\beta}\,) && \Q(\sqrt{\beta'}\,)\\
&\rdLine&\dLine&  \ldLine&\dLine&\rdLine&\dLine &\ldLine \\
& & \Q(\sqrt{D}\,) && F && \Q(\sqrt{d_1}\,) &&\\
&&&\rdLine &  \dLine    &\ldLine&&& \\
&&&& \Q &&&& \\
\end{diagram}

\bigskip

\bigskip

\noindent where $F_c=\Q(\sqrt{d_1},\sqrt{\alpha}\,)$, $F_\gen=\Q(\sqrt{d_1},\sqrt{D}\,)$, and $F=\Q(\sqrt{d_1D}\,)$, for some $\beta,\beta'\in F_c.$

It suffices to show that the rank $r$ of $\Cl_2(\Q(\sqrt{\alpha}\,))$ is $0$.  Since $h_2(D)=1$, we can apply the ambiguous class number formula to the extension $\Q(\sqrt{\alpha}\,)/\Q(\sqrt{D}\,)$:
$$2^r=\#\Am_2(\Q(\sqrt{\alpha}\,)/\Q(\sqrt{D}\,))=\frac{2^{t-1}}{(E_D:H)},$$
where $E_D$ is the unit group of $\Q(\sqrt{D}\,)$,  $H=E_D\cap \mathcal{N}\Q(\sqrt{\alpha}\,)^\times$, with $\mathcal{N}$ denoting the norm from $\Q(\sqrt{\alpha}\,)$ to $\Q(\sqrt{D}\,)$, and $t=t_\fin + t_\infty$, with $t_\fin, t_\infty$, equal to the number of finite, resp.\;infinite primes of $\Q(\sqrt{D}\,)$ ramified in $\Q(\sqrt{\alpha}\,)$.

We claim $t_\fin=1.$  Since in $F_c/F$ there is no ramification outside $\infty$, the only primes ramified in $F_c/\Q$ are $p_1$ and the $p_i|D$, and only with ramification index $2$. Since any prime $p_i|D$ already ramifies in $\Q(\sqrt{D}\,)$, the only prime ideals in $\Q(\sqrt{D}\,)$ ramified in $\Q(\sqrt{\alpha}\,)$ must divide $p_1$. Now $p_1$ splits in $\Q(\sqrt{D}\,)$, say $p_1\OO_{\Q(\sqrt{D}\,)}=\fp_1\fp_1'$. Certainly both of these primes cannot be unramified in $\Q(\sqrt{\alpha}\,)$, for otherwise they are both unramified in $\Q(\sqrt{\alpha'}\,)$. But then also in $F_c$. This implies that $p_1$ is unramified in $F_c/\Q$, a contradiction.  On the other hand, not both $\fp_1$ and $\fp_1'$  are ramified in $\Q(\sqrt{\alpha}\,)$. For otherwise both primes also ramify in $\Q(\sqrt{\alpha'}\,)$ and hence also in $F_\gen$. This implies that the two primes are totally ramified in $F_c/\Q(\sqrt{D}\,)$, contradicting the bound on the ramification index. Therefore, $t_\fin=1$ as claimed.

On the other hand, $t_\infty$ depends on the sign of $\alpha$.

If $\alpha>0$, then $t_\infty=0$; hence $\#\Am_2(\Q(\sqrt{\alpha}\,)/\Q(\sqrt{D}\,))=1$, and so $\Cl_2(\Q(\sqrt{\alpha}\,))$ is trivial and therefore $h_2(\Q(\sqrt{\alpha}\,))=1.$

If $\alpha<0$, then $t_\infty=2$; hence $\#\Am_2(\Q(\sqrt{\alpha}\,)/\Q(\sqrt{D}\,))=4/(E_D:H)$. But recall that $D=d_2$ and so $E_D=\la -1,\varepsilon\ra$, where $\varepsilon$ is the fundamental unit in $\Q(\sqrt{d_2}\,)$, $N\varepsilon=-1$, and $\la \varepsilon^2\ra\subseteq H\subseteq E_D^+=\la \varepsilon^2\ra.$ Therefore $(E_D:H)=4$ and so again $\Cl_2(\Q(\sqrt{\alpha}\,))$ is trivial and therefore $h_2(\Q(\sqrt{\alpha}\,))=1.$

This completes the proof of the lemma.
\end{proof}

By Proposition~\ref{PK4}
$$h_2(M)=\frac1{2}h_2(\Q(\sqrt{\gamma}\,))h_2(\Q(\sqrt{D_1\gamma}\,)).$$
 More explicitly,
$$h_2(M)=\left\{
           \begin{array}{ll}\vspace{.05in}
            \frac1{2}h_2(\Q(\sqrt{d_3d_4\alpha}\,)), & \hbox{if}\;\;D=d_2,\alpha<0, \\\vspace{.05in}
             \frac1{2}h_2(\Q(\sqrt{d_3\alpha}\,))h_2(\Q(\sqrt{d_4\alpha}\,)), & \hbox{if}\;\;D=d_2,\alpha>0, \\
             \frac1{2}h_2(\Q(\sqrt{d_3\alpha}\,))h_2(\Q(\sqrt{d_2d_3\alpha}\,)), & \hbox{if}\;\;D=d_3d_4\,(\Rightarrow \alpha>0).
           \end{array}
         \right.
      $$

\bigskip

\noindent {\bf Assume $D=d_2$.}

\medskip
By the appendices, we see that this  can occur precisely  in the cases
$$(a_1), (b_1),(a_2), (b_2),(a_5), (b_5),(a_6), (b_6),(a_7), (b_7).$$

\begin{prop}\label{P5} Let
\begin{align*}
 \mathcal{A}&=\{(a_1),(a_2),(b_2),(a_5),(b_5),(a_7),(b_7)\}\\
 \mathcal{B}&=\{(b_1),(a_6),(b_6)\}.
\end{align*}
Then
$$h_2(M)\geq\left\{
                                             \begin{array}{ll}
                                               4, & \hbox{if}\;\;\; k\in\mathcal{A}, \\
                                               2, & \hbox{if}\;\;\; k\in\mathcal{B}.
                                             \end{array}
                                           \right.$$
 \end{prop}

\begin{proof} From above, since $D=d_2$, we have
 $$h_2(M)=\left\{
           \begin{array}{ll}\vspace{.05in}
            \frac1{2}h_2(\Q(\sqrt{d_3d_4\alpha}\,)), & \hbox{if}\;\;\alpha<0, \\\vspace{.05in}
             \frac1{2}h_2(\Q(\sqrt{d_3\alpha}\,))h_2(\Q(\sqrt{d_4\alpha}\,)), & \hbox{if}\;\;\alpha>0.
              \end{array}
         \right.
      $$
For both $\alpha>0$ and $\alpha<0$, we can refer to the diagram:

\bigskip

Diagram 3

\begin{diagram}[height=2em,width=2em]
 & & M & &    \\
 &\ldLine &\dLine &\rdLine &    \\
\mathcal{F} & & \Q(\sqrt{\gamma}\,)  & &\Q(\sqrt{D_1\gamma}\,) & &  \\
 &\rdLine & \dLine & \ldLine & \\
 & & \Q(\sqrt{d_2}\,) & &  \\
\end{diagram}

\noindent where $\cF=\Q(\sqrt{d_2},\sqrt{d_3d_4}\,)$ as above.

We will use the ambiguous class number formula on the extensions \\ $\Q(\sqrt{\gamma}\,)/\Q(\sqrt{d_2}\,)$ and $\Q(\sqrt{d_3d_4\gamma}\,)/\Q(\sqrt{d_2}\,)$ to compute the ranks of \\ $\Cl_2(\Q(\sqrt{\gamma}\,))$ and $\Cl_2(\Q(\sqrt{d_3d_4\gamma}\,))$.

\medskip

\noindent {\em Case 1.} Assume $\alpha<0$.

\medskip
  Since $h_2(\Q(\sqrt{\alpha}\,))=1$ and so we only need to consider $r=\rank\Cl_2(\Q(\sqrt{d_3d_4\alpha}\,)).$  Since $h_2(d_2)=1$, we have
$$2^r=\#\Am_2=\#\Am_2\big(\Q(\sqrt{d_3d_4\alpha}\,)/\Q(\sqrt{d_2}\,)\big)=\frac{2^{t-1}}{(E_{d_2}:H)},$$
where $t=t_\fin+t_\infty$. Since $\alpha<0$, \; $t_\infty=2$ and we see as above that $(E_{d_2}:H)=4$.

Now consider $t_\fin$, whose value will depend on the various cases. The rational primes ramified in $\Q(\sqrt{d_3d_4\alpha}\,)/\Q$ are $p_1,p_3,p_4$. By  Appendix~I, we have the follow decompositions of the $p_i$ in $\Q(\sqrt{d_2}\,)$:
\begin{align*}
    p_1\;&\text{splits in all cases,}\\
    p_3\;&\left\{
            \begin{array}{ll}
              \text{splits in cases} & (a_1),(a_5),(b_5),(a_7),(b_7), \\
              \text{is inert in cases} & (b_1),(a_2),(b_2),(a_6),(b_6), \\
             \end{array}\right.\\
    p_4\;&\left\{
            \begin{array}{ll}
              \text{splits in cases} & (a_2),(b_2), \\
              \text{is inert in cases} &  (a_1),(b_1),(a_5),(b_5),(a_6),(b_6),(a_7),(b_7).\\
            \end{array}
          \right.
\end{align*}

Now let $t_\fin=t_1+t_3+t_4,$ where $t_j$ is the number of the prime ideals in $\Q(\sqrt{d_2}\,)$ dividing $p_j$ and ramified in $\Q(\sqrt{d_3d_4\alpha}\,)$. Then by an argument analogous to that in Lemma~\ref{L4}, we have  $t_1=1$. On the other hand, if $p_j$ is inert in $\Q(\sqrt{d_2}\,)$, then (clearly) $t_j=1$. However if for $j=3,4$,  $p_j$ splits, then $t_j=2$: For let $(p_j)=\fp_j\fp_j'$ in $\Q(\sqrt{d_2}\,)$. Then neither $\fp_j$ nor $\fp_j'$ ramifies in $\Q(\sqrt{\alpha}\,)$, cf.\:the proof of Lemma~\ref{L4}. But both ramify in $\mathcal{F}=\Q(\sqrt{d_2},\sqrt{d_3d_4}\,)$, and so each must ramify in $\Q(\sqrt{d_3d_4\alpha}\,)/\Q(\sqrt{d_2}\,)$.
A little calculation then yields
$$t_\fin=\left\{
           \begin{array}{ll}
             4, & \hbox{if}\;\;\;k\in\cA, \\
             3, & \hbox{if}\;\;\;k\in\cB,
           \end{array}
         \right. $$
Hence
$$t=\left\{
           \begin{array}{ll}
             6, & \hbox{if}\;\;\;k\in\cA, \\
             5, & \hbox{if}\;\;\;k\in\cB,
           \end{array}
         \right. $$
which in turn implies that
$$\#\Am_2=\left\{
           \begin{array}{ll}
             2^3, & \hbox{if}\;\;\;k\in\cA, \\
             2^2, & \hbox{if}\;\;\;k\in\cB.
           \end{array}
         \right. $$
Therefore,
$$\rank\Cl_2\big(\Q(\sqrt{d_3d_4\alpha}\,)\big)=\left\{
                                             \begin{array}{ll}
                                               3, & \hbox{if}\;\;\; k\in\mathcal{A}, \\
                                               2, & \hbox{if}\;\;\; k\in\mathcal{B}.
                                             \end{array}
                                           \right.$$
This establishes the proposition for negative $\alpha$.

\bigskip
\noindent {\em Case 2.}  Assume $\alpha>0$.

\medskip
This time we'll compute the ranks $r_j$ of $\Cl_2(\Q(\sqrt{d_j\alpha}\,))$ for $j=3,4.$ We have
$$2^{r_j}=\#\Am_2\big(\Q(\sqrt{d_j\alpha}\,)/\Q(\sqrt{d_2}\,)\big)=\frac{2^{t^{(j)}-1}}{(E_{d_2}:H_j)};$$
as above $(E_{d_2}:H_j)=4$.  Notice that $t^{(j)}=2+t^{(j)}_\fin$ and  $t^{(j)}_\fin=t_{p_1}+t_{p_j}$ ($j=3,4$);
$$t_{p_1}=1\,,\qquad t_{p_j}=\left\{
                           \begin{array}{ll}
                             2, & \hbox{if} \;p_j\;\hbox{splits in}\;\Q(\sqrt{d_2}\,), \\
                             1, & \hbox{if} \;p_j\;\hbox{is inert in}\;\Q(\sqrt{d_2}\,).
                           \end{array}
                         \right. $$
Therefore
\begin{align*}
 p_3\;&\left\{
            \begin{array}{ll}
              \text{splits in cases} & (a_1),(a_5),(b_5),(a_7),(b_7), \\
              \text{is inert in cases} & (b_1),(a_2),(b_2),(a_6),(b_6), \\
             \end{array}\right.\\
    p_4\;&\left\{
            \begin{array}{ll}
              \text{splits in cases} & (a_2),(b_2), \\
              \text{is inert in cases} &  (a_1),(b_1),(a_5),(b_5),(a_6),(b_6),(a_7),(b_7).\\
            \end{array}
          \right.
\end{align*}
We thus have
\begin{align*}
 \rank\Cl_2(\Q(\sqrt{d_3\alpha}\,))=&\left\{
            \begin{array}{ll}
              2\;\text{ in cases} & (a_1),(a_5),(b_5),(a_7),(b_7), \\
              1\;\text{ in cases} & (b_1),(a_2),(b_2),(a_6),(b_6), \\
             \end{array}\right.\\
    \rank\Cl_2(\Q(\sqrt{d_4\alpha}\,))=&\left\{
            \begin{array}{ll}
             2\; \text{ in cases} & (a_2),(b_2), \\
             1\; \text{ in cases} &  (a_1),(b_1),(a_5),(b_5),(a_6),(b_6),(a_7),(b_7).\\
            \end{array}
          \right.
\end{align*}

Considering both these ranks, we see that
$$h_2(M)\geq\left\{
                                             \begin{array}{ll}
                                               4, & \hbox{if}\;\;\; k\in\mathcal{A}, \\
                                               2, & \hbox{if}\;\;\; k\in\mathcal{B}.
                                             \end{array}
                                           \right. $$
 and thus the proposition is established.
\end{proof}

Here are some random examples, using Pari, corroborating  some of the results in Proposition~\ref{P5}:
\begin{align*}
               & \alpha<0\,:\\
    (a_1)\;\; & d_k=8\cdot 17\cdot 47\cdot 3,\;\;\alpha=-5+\sqrt{17},\\
              & \Cl_2(\sqrt{\alpha}\,)\simeq (1),\;\; \Cl_2(\sqrt{47\cdot 3\cdot \alpha})\simeq (2,2,2),\\
    (b_5)\;\; & d_k=61\cdot 5\cdot 11\cdot 4,\;\;\alpha=-9+2\sqrt{5},\\
              & \Cl_2(\sqrt{\alpha}\,)\simeq (1),\;\; \Cl_2(\sqrt{11\alpha})\simeq (2,2,16),\\
    (a_6)\;\; & d_k=8\cdot 17\cdot3\cdot 7,\;\;\alpha=-5+\sqrt{17},\\
              & \Cl_2(\sqrt{\alpha}\,)\simeq (1),\;\; \Cl_2(\sqrt{21 \alpha})\simeq (2,4),\\
    (b_6)\;\; & d_k=41\cdot 5\cdot 7\cdot 4,\;\;\alpha=-11+5\sqrt{5},\\
              & \Cl_2(\sqrt{\alpha}\,)\simeq (1),\;\; \Cl_2(\sqrt{7 \alpha})\simeq (2,8),\\
    (b_1)\;\; & d_k=64\cdot 5\cdot 7\cdot 4,\;\;\alpha=-9+2\sqrt{5},\\
              & \Cl_2(\sqrt{\alpha}\,)\simeq (1),\;\; \Cl_2(\sqrt{7 \alpha})\simeq (2,2);\\
              & \alpha>0\,:\\
    (a_2)\;\; & d_k=8\cdot 113\cdot 3\cdot 7,\;\;\alpha=11+\sqrt{113},\\
              & \Cl_2(\sqrt{-3\alpha}\,)\simeq (2),\;\; \Cl_2(\sqrt{-7\alpha})\simeq (2,2),\\
    (a_5)\;\; & d_k=41\cdot 8\cdot 79\cdot 3,\;\;\alpha=7+2\sqrt{2},\\
              & \Cl_2(\sqrt{-79\alpha}\,)\simeq (2,4),\;\; \Cl_2(\sqrt{-3\alpha})\simeq (2),\\
    (b_5)\;\; & d_k=29\cdot 5\cdot11\cdot 4,\;\;\alpha=7+2\sqrt{5},\\
              & \Cl_2(\sqrt{-11\alpha}\,)\simeq (2,2),\;\; \Cl_2(\sqrt{-\alpha})\simeq (2),\\
    (a_6)\;\; & d_k=8\cdot 41\cdot 3\cdot 7,\;\;\alpha=7+\sqrt{41},\\
              & \Cl_2(\sqrt{-3\alpha}\,)\simeq (2),\;\; \Cl_2(\sqrt{-7\alpha})\simeq (4),\\
    (b_6)\;\; & d_k=401\cdot 5\cdot 3\cdot 4,\;\;\alpha=41+16\sqrt{5},\\
              & \Cl_2(\sqrt{-3\alpha}\,)\simeq (2),\;\; \Cl_2(\sqrt{-\alpha})\simeq (4),\\
    (b_1)\;\; & d_k=29\cdot 5\cdot 3\cdot 4,\;\;\alpha=7+2\sqrt{5},\\
              & \Cl_2(\sqrt{-3\alpha}\,)\simeq (2),\;\; \Cl_2(\sqrt{-\alpha})\simeq (2).
\end{align*}
Observe that  the Main Theorem holds now for all $k$ in cases $\cA$.  But for $k$ in the cases in $\cB$ there's more to do. For  $(a_6),(b_6)$, notice that in the examples that  the product of the $2$-class numbers of  the quartic fields have order $\geq 8$. This is true in general, which follows directly from the next proposition:

\begin{prop}\label{P6}
Assume  and that $k$ satisfies case $(a_6)$ or $(b_6)$. Then
\begin{align*}
h_2\big(\Q(\sqrt{d_3d_4\alpha}\,)\big)&\geq 8\;\;\hbox{if}\;\;\alpha<0,\\
h_2\big(\Q(\sqrt{d_4\alpha}\,)\big)&\geq 4\;\;\hbox{if}\;\;\alpha>0.
\end{align*}
 \end{prop}
\begin{proof} Assume first that $\alpha<0$; by the proof of Proposition~\ref{P5}, $\rank\Cl_2(\Q(\sqrt{d_3d_4\alpha}\,))=2,$ and so there are three unramified quadratic extensions of $$K_0=\Q(\sqrt{d_3d_4\alpha}\,),$$ namely,
$$M=\Q(\sqrt{d_3d_4},\sqrt{\alpha}\,),\;\;M_1=\Q(\sqrt{d_3},\sqrt{d_4\alpha}\,),\;\;M_2=\Q(\sqrt{d_4},\sqrt{d_3\alpha}\,).$$
We claim that  at least one (hence all) of these three fields has odd class number equal to half that of $K_0$, and hence that the $2$-class field tower of $K_0$ terminates at $K_0^1$, see Proposition~7 in \cite{BLS}. For instance, applying Kuroda's class number formula to the $V_4$-extension $M/\Q(\sqrt{d_2}\,)$ yields
$$\frac{h_2(M)}{h_2(K_0)}=\frac1{2}q(M/\Q(\sqrt{d_2}\,) h_2(\mathcal{F})h_2(\Q(\sqrt{\alpha}\,))/h_2(d_2)^2=\frac1{2}q(M/\Q(d_2),$$
since $h_2(d_2)=h_2(\Q(\sqrt{\alpha}\,))=h_2(\mathcal{F})=1.$ Moreover, since $M/\Q(\sqrt{d_2}\,)$ is  a $V_4$-extension of CM-fields, we see that
$$q=q(M/\Q(\sqrt{d_2}\,)=\frac{Q(M)}{Q(K_0)Q(\Q(\sqrt{\alpha}\,))}=1$$
 (recall that $Q(M)=1$ and of course $q\in\Z$). Thus $h_2(M)=\frac1{2}h_2(K_0),$ which further  implies $(i=1,2)$
$$h_2(M_i)=\frac1{2}h_2(K_0),$$
as desired.

Now we claim that $\rank\Cl_2(M_2)=2$. This will establish the proposition, since then $h_2(M_2)\geq 4$ implying that $h_2(K_0)\geq 8$. For starters consider the following diagram:

\bigskip
Diagram 4

\begin{diagram}[height=2em,width=2em]
 & & M_2 & &    \\
 &\ldLine &\dLine &\rdLine &    \\
\mathcal{F}_1 & & K_0  & & K_1 & &  \\
 &\rdLine & \dLine & \ldLine & \\
 & & \Q(\sqrt{d_2}\,) & &  \\
\end{diagram}

\noindent where $\cF_1=\Q(\sqrt{d_2},\sqrt{d_4}\,)$ and $K_1=\Q(\sqrt{d_3\alpha}\,)$. First notice that $h_2(\cF_1)=1$, for by the ambiguous  class number formula
$$2^{\rank\Cl_2(\cF_1}=\#\Am_2(\cF_1/\Q(\sqrt{d_2}\,))=2^{t-1}/(\la -1,\varepsilon\ra:H)=1,$$
since $t=t_\infty+t_\fin$ with $t_\infty=2$, $t_\fin=1$, and $(\la -1,\varepsilon\ra:H)=2^2$, as is easily verified. Hence we can use the ambiguous class number formula to calculate the rank $r$ of $\Cl_2(M_2)$:
$$2^r=\#\Am_2(M_2/\cF_1)=\frac{2^{t-1}}{(E_{\cF_1}:H)},$$
where $t=t_\fin+t_\infty=t_\fin$, since $t_\infty=0$ ($\cF_1$ is already totally complex). We see that $t_\fin=t_1+t_3$, where $t_j$ is the number of prime ideals in $\cF_1$ dividing $p_j$ and ramifying in $M_2$. Now $p_1$ splits completely in $\cF_1$ and hence there are two prime ideals in $\cF_1$ above the unique prime ideal in $\Q(\sqrt{d_2}\,)$ that ramifies in $K_0$. Therefore $t_1=2$. On the other hand, $p_3$ splits into two prime ideals in $\cF_1$, each of which is ramified in $M_2$. Therefore $t_3=2.$ Thus $t=4$ and so
$$\#\Am_2(M_2/\cF_1)=\frac{2^3}{(E_{\cF_1}:H)}.$$

Now we claim that $(E_{\cF_1}:H)=2$, where
$$E_{\cF_1}=\left\{
              \begin{array}{ll}
                \la -1,\varepsilon_2, E_{\cF_1}^2\ra, & \hbox{in the case}\;\;(a_6) \\
                \la i,\varepsilon_2\ra, & \hbox{in the case}\;\;(b_6)
              \end{array}
            \right.$$
where $\varepsilon_2$ is the fundamental unit of $\Q(\sqrt{d_2}\,)$
and $H=E_{\cF_1}\cap\cN$, with $\cN=N_{M_2/\cF_1}M_2^\times.$ It is easy to see that $(E_{\cF_1}:H)\leq 4$. We'll  show that the group of roots of unity $W_{\cF_1}\subseteq H$, which will yield $(E_{\cF_1}:H)\leq 2$. Toward this end, we  use Hasse's local norm theorem: $u\in H$ if and only if $u\in\cN_{\fP}$ for all the finite and infinite primes $\fP$ in $M_2$, where $\cN_{\fP}= N_{(M_2)_\fP/(\cF_1)_\fp} (M_2)_\fP^\times$, and  where $\fp=\fP\cap \cF_1$,    for relevant completions. Since containment holds for the unramified primes, we need only check the two primes $\fp_1, \fp_1'$ in $\cF_1$ above $p_1$ that ramify in $M_2$, and similarly for $\fp_3,\fp_3'$ above $p_3$.

We'll start with $\fp=\fp_1$ or $\fp'$ above $p_1$. First suppose $p_1\not=2$. Then
$$u\in \cN_{\fP}\;\;\Leftrightarrow\;\;u\equiv \xi^2\bmod \fp\OO_{(\cF_1)_\fp},\;\; \text{for some}\;\xi\in \OO_{(\cF_1)_\fp},$$
see the proof of Proposition~15 in \cite{BS1}.
Now, in case $(a_6)$ \;$-1\equiv a^2\bmod p_1$, since $p_1\equiv 1\bmod 4$ and thus $-1\in\cN_{\fP}$. In case $(b_6)$, $i\equiv \xi^2\bmod \fp\OO_{(\cF_1)_\fp}$, since $-1\equiv a^4\bmod p_1$, because $p_1\equiv 1\bmod 8$. Therefore $W_{(\cF_1)_\fp}\subseteq \cN_{\fP}$ for all $\fp$ dividing $p_1$.

If $p_1=2$\;(only in case $(a_6)$), then $u\in\cN_\fP$ if and only if $u^2\equiv 1\bmod 16\OO_{(\cF_1)_{\fp}}$. Thus for $u=-1$, the congruence is trivially true; so again $W_{(\cF_1)_\fp}\subseteq \cN_{\fP}$ for all $\fp$ dividing $p_1$.

Next, consider $\fp\mid p_3.$ If $p_3\not=2$, then since $p_3$ is inert in $F_2=\Q(\sqrt{d_2}\,)$ we see that $\#\OO_{F_2}/p_3\OO_{F_2}=p_3^2\equiv 1\bmod 8$. Therefore $-1\equiv a^4\bmod p_3\OO_{F_2},$ and so $-1\in\cN$ and $i\in\cN$ for case $(b_6)$. For $p_3=2,$ the argument above shows that $-1\in\cN$.

Therefore $W_{\cF_1}\subseteq H,$ and so $(E_{\cF_1}:H)\leq 2$. Thus $\rank\Cl_2(M_2)\geq 2,$ and hence $=2$, since $\rank\Cl_2(K_0)=2$ (recall that the $2$-class field tower of $K_0$ terminates at $K_0^1$). Thus the proposition is true for $\alpha<0$.

\medskip
Assume $\alpha>0$. We  consider  Diagram 4, where  $\cF_1=\Q(\sqrt{d_2},\sqrt{d_4}\,)$ as above, but this time
$$M_2=\Q(\sqrt{d_4},\sqrt{\alpha}\,),\quad K_0=\Q(\sqrt{d_4\alpha}\,), \quad K_1=\Q(\sqrt{\alpha}\,).$$
Since $M_2/K_0$ is unramified, if we can show that $h_2(M_2)\not=1$, then we'll be done. We'll do this by showing that $r=\rank\Cl_2(M_2)=1.$ By the ambiguous class number formula (recall that $h_2(\cF_1)=1$)
$$2^r=\#\Am_2(M_2/\cF_1)=\frac{2^{t-1}}{(E_{\cF_1}:H)},$$
where $t=t_\fin$, since $t_\infty=0$. This time only the two primes  in $\cF_1$ above the prime $p_1$ in $\Q(\sqrt{d_2}\,)$ ramified in $K_0$ and $K_1$ ramify in $M_2/\cF_1$ and so we have $t_\fin=t_{1}=2.$

On the other hand, as in the argument above, we see that $W_{\cF_1}\subseteq H$. Now, we claim that $\varepsilon_2\in H$. For if $p_1$ is odd, then we need to show that
$$\varepsilon_2\equiv a^2\bmod \fp_1,$$
for $\fp_1\in \cF_1$, dividing $p_1$. This congruence holds if $(p_1/p_2)_4(p_2/p_1)_4=+1.$ But this is true since $h_2(d_1d_2)\geq 4$, see, for example p.\,56 in \cite{Lem}. Hence $\varepsilon_2\in H$.  A similar argument works for $p_1=2.$

Thus  $\#\Am_2(M_2)=2$, since $t=2$ and $(E_{\cF_1}:H)=1$. Therefore \\ $\rank\Cl_2(M_2)=1$, implying that $h_2(\sqrt{d_4\alpha}\,)\geq 4.$

This establishes the proposition.
\end{proof}

Notice that since $h_2(\Q(\sqrt{d_3\alpha}\,)\not=1$, \;$h_2(M)\geq 4$.

We now have verified the Main Theorem for all $k$ for which $D=d_2$ and $k$ in (A) (with the usual exception of case $(b_1)$).

\bigskip

\noindent {\bf Assume $D=d_3d_4$}.

By the appendices this occurs for the cases $(a_3),(b_3),(a_4),(b_4),(a_8),(b_8)$.

\begin{prop}\label{P7} For all $k$ for which $D=d_3d_4$,\; $h_2(M)\geq 4.$
\end{prop}
\begin{proof} Consider the diagram

\bigskip
Diagram 5

\begin{diagram}[height=2em,width=2em]
 & & M & &    \\
 &\ldLine &\dLine &\rdLine &    \\
\mathcal{F} & & \Q(\sqrt{d_3\alpha}\,)  & &\Q(\sqrt{d_2d_3\alpha}\,) & &  \\
 &\rdLine & \dLine & \ldLine & \\
 & & \Q(\sqrt{d_3d_4}\,) & &  \\
\end{diagram}

\noindent where once again $\cF=\Q(\sqrt{d_2},\sqrt{d_3d_4}\,)$ and $M=\Q(\sqrt{d_3\alpha},\sqrt{d_2d_3\alpha}\,)$.

It suffices to show that
$$r=\rank\Cl_2(\Q(\sqrt{d_3\alpha}\,))=1\,,\;\;\rank\Cl_2(\Q(\sqrt{d_2d_3\alpha}\,))\geq 2.$$
As usual, we'll apply the ambiguous class number formula (notice that $h_2(d_3d_4)=1$):
$$2^r=\#\Am_2=\#\Am_2(\Q(\sqrt{d_3\alpha}\,)/\Q(\sqrt{d_3d_4}\,))=\frac{2^{t-1}}{(E:H)}\,,$$
where $E=E_{\Q(\sqrt{d_3d_4}\,)}=\la -1,\varepsilon\ra$, for $\varepsilon$  the fundamental unit of $\Q(\sqrt{d_3d_4}\,)$. We have $t=t_\fin+t_\infty$, with $t_\infty=2.$ Hence $E^2\subseteq H\subseteq E^+\subseteq E$, where $E^+$ is the subgroup of totally positive units in $E$ and thus $=\la \varepsilon\ra$ since $N\varepsilon=1$. This implies that $(E:E^+)=2$ which in turn implies that $(E:H)=2$ or $4$, according as $\varepsilon$ is or is not, resp.,  a norm from $\Q(\sqrt{d_3\alpha}\,)$ to $\Q(\sqrt{d_3d_4}\,)$. (We'll see that $(E:H)=2.$)

Now we consider $t_\fin$. The only primes of $\Q(\sqrt{d_3d_4}\,)$ ramified in $\Q(\sqrt{d_3\alpha}\,)$ are those dividing $p_1$ and hence, once again, $t_{p_1}=1.$ Thus $t=3$ and so $\#\Am_2=2^2/(E:H)\leq 2$. On the other hand, $\Q(\sqrt{d_3},\sqrt{\alpha}\,)$ is an unramified quadratic extension of $\Q(\sqrt{d_3\alpha}\,)$ and hence $h_2(\Q(\sqrt{d_3\alpha}\,))\geq 2$. This implies that $\#\Am_2>1$ and thus $=2$. Therefore,
$\rank\Cl_2(\Q(\sqrt{d_3\alpha}\,))=1\,.$

Now consider $\Cl_2(\Q(\sqrt{d_2d_3\alpha}\,))$. This time
$$\#\Am_2(\Q(\sqrt{d_2d_3\alpha}\,)/\Q(\sqrt{d_3d_4}\,))=\frac{2^{t-1}}{(E:H)}\,,$$
with $t=t_\fin+2$ and as above $(E:H)=2$ or $4$. The primes of $\Q(\sqrt{d_3d_4}\,)$ that can possibly ramify in $\Q(\sqrt{d_2d_3\alpha}\,)$ are those above $p_1$ and $p_2$. As before, $t_{p_1}=1.$ On the other hand,

\begin{align*}
 p_2\;&\left\{
            \begin{array}{ll}
              \text{splits in cases} & (a_8),(b_8), \\
              \text{is inert in cases} & (a_3),(b_3),(a_4),(b_4). \\
             \end{array}\right.
\end{align*}
By the diagram, if $\fp$ is a prime in $\Q(\sqrt{d_3d_4}\,)$ dividing $p_2$, then $\fp$ ramifies in $\cF$, but not in $\Q(\sqrt{d_3\alpha}\,).$ Hence $\fp$ ramifies in $\Q(\sqrt{d_2d_3\alpha}\,)$. Consequently,

\begin{align*}
 t_{p_2}=&\left\{
            \begin{array}{ll}
             2\;\; \text{ in cases} & (a_8),(b_8), \\
             1\;\; \;\text{in cases} & (a_3),(b_3),(a_4),(b_4). \\
             \end{array}\right.
\end{align*}
This implies that
$$\#\Am_2(\Q(\sqrt{d_2d_3\alpha}\,)/\Q(\sqrt{d_3d_4}\,))=\left\{
                                                           \begin{array}{ll}
                                                             2^4/(E:H), & \hbox{in cases}\;\; (a_8),(b_8), \\
                                                             2^3/(E:H), & \hbox{in cases}\;\; (a_3),(b_3),(a_4),(b_4).
                                                           \end{array}
                                                         \right.$$
In cases $(a_8),(b_8)$, $\#\Am_2\geq 2^2$ and thus $\rank\Cl_2(\Q(\sqrt{d_2d_3\alpha}\,))\geq 2$, as desired.

Now consider the cases $(a_3),(b_3),(a_4),(b_4)$. In these cases we claim that\\ $(E:H)=2.$ To see this, we will show that $\varepsilon\in N_{K/F}K^\times$ where for the moment we are letting $K=\Q(\sqrt{d_2d_3\alpha}\,)$ and $F=\Q(\sqrt{d_3d_4}\,)$. By Hasse's local norm theorem $\varepsilon\in N_{K/F}K^\times$  if and only if $\varepsilon\in N_{K_\fP/F_\fp}K_\fP^\times$ for all primes $\fP$ (finite and infinite) in $K$ where $\fp=\fP\cap F$. For $\fP/\fp$ unramified we have $\varepsilon\in N_{K_\fP/F_\fp}K_\fP^\times$. For all the infinite primes this containment holds since $\varepsilon$ is totally positive in $F$.  Hence we need only consider the finite primes $\fp_1$, one of the two primes above $p_1$ in $F$, and $\fp_2$, the inert prime above $p_2$. Suppose $p_1$ is odd; then for $\fP_1$ in $K$ dividing $\fp_1$,
$\varepsilon\in N_{K_{\fP_1}/F_{\fp_1}}K_{\fP_1}^\times$ if and only if $\varepsilon\equiv \xi^2\bmod \fp_1\OO_{F_{\fp_1}}.$
But this congruence already holds for the case $\Q(\sqrt{d_3\alpha}\,)/F$ from the previous case above. A similar argument works for $p_1=2$. Hence $\varepsilon$ is a norm from $K_{\fP_1}$ to $F_{\fp_1}$. By Hilbert's product formula we thus have
$\varepsilon\in N_{K_{\fP_2}/F_{\fp_2}}K_{\fP_2}^\times$. Therefore $\varepsilon\in N_{K/F}K^\times.$ This implies that  $(E:H)=2$ whence $\#\Am_2=2^2.$  Therefore $\rank\Cl_2(\Q(\sqrt{d_2d_3\alpha}\,))=2$.

All this establishes the proposition.

\end{proof}

\bigskip

Here are some examples  corroborating  the results in this proposition:
\begin{align*}
    (a_3)\;\; & d_k=8\cdot 5\cdot 7\cdot 79,\;\;\alpha=(25+\sqrt{553})/2,\\
              & \Cl_2(\sqrt{-7\alpha}\,)\simeq (8),\;\; \Cl_2(\sqrt{-35\alpha})\simeq (2,16),\\
    (a_4)\;\; & d_k=8\cdot 5\cdot 31\cdot 7,\;\;\alpha=15+\sqrt{217},\\
              & \Cl_2(\sqrt{-31\alpha}\,)\simeq (2),\;\; \Cl_2(\sqrt{-31\cdot 5\cdot \alpha})\simeq (2,2),\\
    (a_8)\;\; & d_k=5\cdot 17\cdot 11\cdot 31,\;\;\alpha=(19+\sqrt{341})/2,\\
              & \Cl_2(\sqrt{-11\alpha}\,)\simeq (2),\;\; \Cl_2(\sqrt{-11\cdot 17\cdot \alpha})\simeq (2,2),\\
    (b_8)\;\; & d_k=17\cdot 5\cdot 19\cdot 4,\;\;\alpha=6+\sqrt{19},\\
              & \Cl_2(\sqrt{-19\alpha}\,)\simeq (4),\;\; \Cl_2(\sqrt{-19\cdot 5\cdot \alpha})\simeq (2,4).\\
\end{align*}

\bigskip

Summarizing the results thus far:

{\em The Main Theorem is valid for all $k$ satisfying the cases $(a_i)$ for $i=1,\dots, 8$, and $(b_i)$ for $i=2,\dots,8$.}

\bigskip
\section{{\bf (B) the cases for which $G$ is quaternion of order $8$}}

By Appendix~I,\; $G\simeq Q$ precisely in cases $(a_9),(b_9), (a_{10}),(b_{10}).$ But since for $(b_9)$, $G^+/G_3^+\simeq 32.037$, we know that the length of the narrow $2$-class field tower of $k$ is $2$. Hence we will only consider the cases
$(a_9),(a_{10}),(b_{10}).$ In these three cases $d_k=d_1\cdot d_2\cdot d_3d_4$ is an $H_8$-factorization, i.e.,
$$(d_1d_2/p_3)=(d_1d_2/p_4)=(d_2d_3d_4/p_1)=(d_3d_4d_1/p_2)=+1,$$
and $k^2=k^1(\sqrt{\mu}\,)$, where $\mu\in\cF=\Q(\sqrt{d_1},\sqrt{d_2}\,)$ with $\mu=\beta u_2$ such that $\beta=x_1\sqrt{d_1}+x_2\sqrt{d_2}\in \OO_\cF$ not divisible by any rational primes, and $u_2\in\Z[\sqrt{d_2}]$ is a unit of norm $-1$, and where $x_1,x_2,x_3$ are rationals satisfying $x_1^2d_1-x_2^2d_2=-d_3d_4x_3^2$, cf.\,\S4 in \cite{Lem}, or \cite{Lem4}. (Slightly more detail: Without loss of generality let $d_2=p_2\equiv 1\bmod 4$. Then $p_2=s^2+t^2$,  with $t$ odd, and by \cite{Lem4} $\mu=\beta(s+\sqrt{p_2}),$ with $s^2-p_2=-t^2.$ But then the ideal $\fa=(s+\sqrt{p_2},t)$ satisfies $\fa^2=(s+\sqrt{p_2})$, as is easily seen. Hence $\fa$ is principal, say $=(\alpha)$ in $\Q(\sqrt{p_2}\,)$, since $h_2(p_2)=1$. Thus $(s+\sqrt{p_2})=(\alpha^2),$ and so $s+\sqrt{p_2}=\eta \alpha^2$ for some unit $\eta$ in $\Q(\sqrt{p_2}\,)$. Moreover, $N\eta=N(s+\sqrt{p_2})/(N\alpha)^2=-(t/N\alpha)^2<0,$ implying that $N\eta=-1$. Take $u_2=\eta$, if $\eta\in\Z[\sqrt{p_2}]$, otherwise let $u_2=\eta^3.$)  Notice  that $\mu$ is totally positive in $k^1$.

Analogous to the previous section, let $L=k^1(\sqrt{d_3},\sqrt{\mu}\,)$. Then we have the following $V_4$-extension of CM-fields:
$$L=k^1(\sqrt{d_3},\sqrt{\mu}\,)$$
$$k_+^1=k^1(\sqrt{d_3}\,),\;\;k^2=k^1(\sqrt{\mu}\,),\;\;N=k^1(\sqrt{d_3\mu}\,)$$
$$k^1=\Q(\sqrt{d_1},\sqrt{d_2},\sqrt{d_3d_4}\,).$$

We need only show that $h_2(L)/h_2(k_+^1)\geq 1.$ By Kuroda's class number formula, we have
$$\frac{h_2(L)}{h_2(k_+^1)}=\frac1{2^{1+v}}q(L/k^1)h_2(k^2)h_2(N)/h_2(k^1)^2\,.$$
Now,
$$\frac1{2^{1+v}}q(L/k^1)= \frac1{2}\frac{Q(L)}{Q(k_+^1)Q(N)},\;\; h_2(k^2)=1,\;\;h_2(k^1)=2\,,$$
and therefore,
$$\frac{h_2(L)}{h_2(k_+^1)}=\frac1{8}\frac{Q(L)}{Q(k_+^1)Q(N)}\, h_2(N)\,.$$
By arguments similar to those in Lemma~\ref{L3.2}, we have $Q(L)=Q(k_+^1)=2$. Since $Q(N)\leq 2$, we will be done if we show that $h_2(N)\geq 16.$ (With a little more effort, we could show that $Q(N)=1$, but this is not needed here.) Toward this end, we start with a lemma.
\begin{lem}\label{L5} Let $\cF=\Q(\sqrt{d_1},\sqrt{d_2}\,).$ Then
$$h_2(N)\geq h_2(\cF(\sqrt{d_3\mu}\,))h_2(\cF(\sqrt{d_4\mu}\,)).$$
\end{lem}
\begin{proof} Consider the following $V_4$-extension of CM-fields:
$$N=k^1(\sqrt{d_3\mu}\,)=k^1(\sqrt{d_4\mu}\,)$$
$$k^1=\cF(\sqrt{d_3d_4}\,),\;\;\cF_3=\cF(\sqrt{d_3\mu}\,),\;\; \cF_4=\cF(\sqrt{d_4\mu}\,)$$
$$\cF=\Q(\sqrt{d_1},\sqrt{d_2}\,).$$
Hence, by Kuroda's class number formula,
$$\frac{h_2(N)}{h_2(k^1)}=\frac1{2^{1+v}}q(N/\cF)h_2(\cF_3)h_2(\cF_4)/h_2(\cF)^2.$$
By Kuroda's class number formula applied to the extension $\cF/\Q$, it can be seen that $h_2(\cF)=1$. Moreover, we also have
$h_2(k^1)=2$ and \\$q(N/\cF)/2^{1+v}=1/2\cdot Q(N)/(Q(\cF_3)Q(\cF_4)).$ But $Q(\cF_i)=1$ for $i=3,4$ since $\cF_i/\cF$ are essentially ramified. Therefore,
$$h_2(N)=Q(N)h_2(\cF_3)h_2(\cF_4)\geq h_2(\cF_3)h_2(\cF_4).$$
\end{proof}
\begin{prop}\label{P8} For $\cF_3$ and $\cF_4$ as above,
$$\rank\Cl_2(\cF_3)=\rank\Cl_2(\cF_4)\geq 2.$$
Therefore $h_2(N)\geq 16$.
\end{prop}
\begin{proof} Since $\cF$ has odd class number, we may use the ambiguous class number formula on $\cF_i/\cF$ for $i=3,4$. We have
$$\#\Am_2(\cF_i/\cF)=\frac{2^{t-1}}{(E_\cF:H)},$$
where $t=t_\fin+t_\infty$. Since $\cF$ is totally real and $\cF_i$ are totally complex, we see that $t_\infty=4$. Also $(E_\cF:H)\leq (E_\cF:E_\cF^2)=2^4.$  Now we consider $t_\fin$. But first we claim that $\cF_3$ and $\cF_4$ are conjugate over $\Q$: For let $\Gal(\cF/\Q)=\la \sigma_1,\sigma_2\ra$ where $\sigma_i(\sqrt{d_j})=\left\{
                                                    \begin{array}{ll}
                                                      -\sqrt{d_j}, & \hbox{if}\;\; j=i, \\
                                                      \;\;\,\sqrt{d_j}, & \hbox{if}\;\; j\not=i.
                                                    \end{array}
                                                  \right.$
We now have
$$
\mu^{1+\sigma_1}=d_3d_4x_3^2u_2^2,\qquad
\mu^{1+\sigma_2}=d_3d_4x_3^2,
$$
whence
$$\sqrt{\mu^{\sigma_1}}=\frac{\sqrt{d_3d_4}x_3u_2}{\sqrt{\mu}},\qquad \sqrt{\mu^{\sigma_2}}=\frac{\sqrt{d_3d_4}x_3}{\sqrt{\mu}}.$$
This implies that $\cF_3^{\sigma_i}=\cF_4.$ Hence $\Cl_2(\cF_3)\simeq \Cl_2(\cF_4)$.

Now consider the diagram

\newpage

Diagram 6

\begin{diagram}[height=2em,width=2em]
   &    &     & \cK &   &   &  &   &  & \\
 & & \ldLine(3,2)   & \dLine & \rdLine(3,2)     &   &     &   & & \\
 \cF_3 & &  & \cF_4 &  &  & k^1  & & &\\
 & \rdLine(3,2)& & \dLine & & \ldLine(3,2) & \dLine  & & & \\
  &  & &\cF & &  & &  & & \\
 &  & & \dLine & &  &   &  &  & \\
   & &  & &  &  & k & & &\\
 &   &  & & & \ruLine(3,2)        &     &   & & \\
   &    &     & \Q &   &   &  &   &  & \\
\end{diagram}

\noindent where $\cK=k^1(\sqrt{d_3\mu}\,)=k^1(\sqrt{d_4\mu}\,)$, $\cF_j=\cF(\sqrt{d_j\mu}\,)$ and $\cF=\Q(\sqrt{d_1},\sqrt{d_2}\,).$  We  have by Appendix~I
\begin{align*}
 p_3 \;&\left\{
            \begin{array}{ll}
              \text{splits completely in $\cF$ for case} & (a_9), \\
              \text{splits into two primes in $\cF$ for cases} &  (a_{10}),(b_{10}) \\
             \end{array}\right.\\
  p_4\;&\left\{
            \begin{array}{ll}
              \text{splits into two primes in $\cF$ for case} & (a_9), \\
              \text{splits completely in $\cF$ for cases} &  (a_{10}),(b_{10}).\\
            \end{array}
          \right.
\end{align*}

Consider the case $(a_9)$; then $p_3\OO_\cF=\fp_1\fp_2\fp_3\fp_4$ for distinct prime in $\cF$. If no $\fp_i$ ramifies in $\cF_3$, hence also in $\cF_4$, then $p_3$ would be unramified in $\cK$ and thus also in $k$, which can't happen. On the other hand, not all the $\fp_i$ can ramify in $\cF_3$, otherwise they also ramify in $\cF_4$ as well as $k^1$. This implies they are totally ramified in $\cK$, with ramification index $4$; but $p_3$ only ramifies in $k$, hence with ramification index $2$, a contradiction. Hence two of the $\fp_i$ ramify in $\cF_3$; the other two in $\cF_4$. Thus $t_{p_3}=2$. In a similar manner, $t_{p_4}=1.$ Therefore $t_\fin=3$ and so $t=7$. Hence
$$\#\Am_2(\cF_i/\cF)=\frac{2^6}{(E_\cF:H)}\geq 2^2,$$
whence $\rank\Cl(\cF_i)\geq 2.$

The cases $(a_{10}), (b_{10})$ are handled in a completely analogous way and are left to the reader.

From this, we see that $h_2(\cF_i)\geq 4$; therefore, $h_2(N)\geq 16$, as desired.
\end{proof}

\bigskip
Here is an example of a field in case $(a_9)$:
$$d_k=13\cdot 5\cdot 131\cdot 7,\quad \beta=4\sqrt{13}+15\sqrt{5},\quad u_2=2+\sqrt{5},$$
$$\Cl_2(N)\simeq (2,2,4,4),\quad \Cl_2(k_+^1)\simeq (2,2,2,2,8).$$

\bigskip

Summarizing the results thus far:

{\em The Main Theorem is valid for all $k$ satisfying the cases $(a_i)$ for $i=1,\dots, 10$, and $(b_i)$ for $i=2,\dots,8$, and $10$.}

\bigskip
\section{{\bf (C) the cases for which $G$ is abelian, except for $G^+/G^+_3\simeq 64.150$}}

These are all the cases $(c_i)$ and $(d_i)$ except $(c_3), (d_1)$. For many of these cases we know that $k$ has narrow $2$-class field tower of length $2$. This leaves  $(c_2),(d_5),$ and $(d_8)$ to consider. Notice that these are the cases where $G^+/G_3^+\simeq 32.033$. We will show that for these  the tower length is at least $3$.

Now consider $(c_2),(d_5),(d_8)$. We'll start with $(c_2)$. Let $F=\Q(\sqrt{d_1d_2d_3}\,)$; then $d_F=d_2\cdot d_1d_3$ is a $C_4$-splitting of $d_F$. This implies that there is an unramified cyclic quartic extension $\cF$ of $F$. Thus there is a primitive rational solution $(a,b,c)$ satisfying $a^2=b^2d_2+c^2d_1d_3$ such that $\alpha=a+c\sqrt{d_1d_3}$ with $\cF=F(\sqrt{\alpha}\,)$. Moreover, $\alpha\alpha'=a^2-c^2d_1d_3=b^2d_2<0$. Hence without loss of generality choose $\alpha>0$ and so $\alpha'<0$.

Now consider the extension $\cF/\Q(\sqrt{d_1d_3}\,).$ Then by shifting by \\ $k^1=\Q(\sqrt{d_1d_4},\sqrt{d_2d_4},\sqrt{d_3d_4}\,)$ we obtain the finitely unramified $V_4$-extension $L/k^1$ with $L=k^1\cF=k^1(\sqrt{d_2},\sqrt{\alpha}\,)$:
$$L=k^1(\sqrt{d_2},\sqrt{\alpha}\,)$$
$$k_+^1=k^1(\sqrt{d_2}\,),\qquad K=k^1(\sqrt{\alpha}\,),\quad K'=k^1(\sqrt{\alpha'}\,)=k^1(\sqrt{d_2\alpha}\,)$$
$$k^1$$
Applying Kuroda's class number formula for $V_4$-extensions to $L/k^1$ yields
$$\frac{h_2(L)}{h_2(k_+^1)}=\frac1{2^{1+v}}q(L/k^1)h_2(K)^2,$$
where $q(L/k^1)=(E_{L}:E_{k_+^1}E_{K}E_{K'})$ and $v=0$ or $1$ with $v=0$ unless $L=k^1(\sqrt{\eta_1},\sqrt{\eta_2}\,)$ for some $\eta_j\in E_{k^1}$, and where we recall that $h_2(k^1)=1$.

We want to show that $h_2(L)/h_2(k_+^1)\geq 1$. Notice that it suffices to show that $h_2(K)\geq 2.$  To this end we apply the ambiguous class number formula to the extension $K/k^1$, observing once again that $h_2(k^1)=1$;
$$2^{\rank\Cl_2(K)}=\#\Am_2(K/k^1)=\frac{2^{t-1}}{(E_{k^1}:H_{K/k^1})},$$
where $t=t_\fin+t_\infty$ is the sum of the number of finite, resp. infinite, places of $k^1$  that ramify in $K$, and $H=H_{K/k^1}=E_{k^1}\cap N_{K/k^1}K^\times$. Since $K/k^1$ is finitely unramified, $t_\fin=0$, and hence we only need to consider the infinite places of $k^1$ that ramify in $K$. We can see that  $t=t_\infty=\#\{\sigma\in\Gal(k^1/\Q):\alpha^\sigma<0\}$. Now $\Gal(k^1/\Q)=\la \sigma_1,\sigma_2,\sigma_3\ra$ where
$$\sigma_i(\sqrt{d_jd_4}\,)=\left\{
                              \begin{array}{ll}
                                -\sqrt{d_jd_4}, & \; \hbox{if}\;\; j=i, \\
                                \phantom{-}\sqrt{d_jd_4}, & \; \hbox{if}\;\; j\not=i.
                              \end{array}
                            \right.$$
Since $\Q(\alpha)=\Q(\sqrt{d_1d_3}\,)$ we see that
$$\alpha^\sigma <0\;\;\text{iff}\;\;\sigma\in S=\{\sigma_1,\sigma_3, \sigma_1\sigma_2,\sigma_2\sigma_3\}.$$
Thus $t=4$ and so $2^{t-1}=2^3.$

Now consider $H.$ By local-global properties, we have
$$H=\{\varepsilon\in E_{k^1}:\varepsilon^\sigma>0, \;\sigma\in S\}.$$

We now need to compute $E_{k^1}$ in order to determine the signs of the conjugates of the units. For starters, since $k^1$ is a multi-quadratic extension of $\Q$ of degree $8$ and all the quadratic subfields of $k^1$ have odd class number except for $k$, for which $h_2(k)=4$, we have (see, for example, p.\,43 in \cite{BLS2})
$$1=h_2(k^1)=2^{-9}q(k^1/\Q)\prod_{{F\subseteq k^1,\,}\\{[F:\Q]=2}}h_2(F)=2^{-7}q(k^1/\Q),$$
which implies that $q(k^1/\Q)=2^7$, where $q(k^1/\Q)=(E_{k^1}:E)$, with $E$ the  subgroup of $E_{k^1}$ generated by all the units in the quadratic subfields of $k^1$, i.e.,
$$E=\la -1,\varepsilon_{12},\varepsilon_{13},\varepsilon_{14},\varepsilon_{23},\varepsilon_{24},\varepsilon_{34},\varepsilon_k \ra,$$
where $\varepsilon_{ji}=\varepsilon_{ij}$ is the fundamental unit in $\Q(\sqrt{d_id_j}\,)$, etc.  Now we have the following useful lemmas. In the first lemma $\varepsilon_m$ denotes the fundamental unit of $\Q(\sqrt{m}\,)$.

\begin{lem}\label{L6} Let $p,q$ be distinct primes $\equiv 3\bmod 4$. Then

\medskip
\noindent i) $\sqrt{\varepsilon_{pq}}=\frac1{2}(a\sqrt{p}+b\sqrt{q})$,\; for some  $a,b\in\Z$ with $a\equiv b\bmod 2$;

moreover, $\frac1{4}(a^2p-b^2q)=(p/q)=-(q/p)$;

\medskip
\noindent ii) $\sqrt{\varepsilon_{2p}}=a\sqrt{2}+b\sqrt{p}$,\; for some  $a,b\in\Z$;

moreover, $a^22-b^2p=(2/p)=(-p/2)$;

\medskip
\noindent iii) $\sqrt{\varepsilon_{p}}=\frac1{2}(a\sqrt{2p}+b\sqrt{2})$,\; for some odd $a,b\in\Z$;

moreover, $\frac1{2}(a^2p-b^2)=-(-p/2)$.

\end{lem}
\begin{proof}
To see this (see for example Proposition~3 in \cite{BLS},  or see  \cite{Kub}),  since  $\varepsilon=\varepsilon_{pq}, \varepsilon_{2p},$ and $\varepsilon_{p}$  have positive norm, we have (as a mini-version of Hilbert's Theorem 90) $\varepsilon=(1+\varepsilon)/(1+\varepsilon')$, thus $\varepsilon=(1+\varepsilon)^2/N(1+\varepsilon)$, which implies that  $\Q(\sqrt{\varepsilon}\,)=\Q(\varepsilon,\sqrt{N(1+\varepsilon)}\,)$; recall that the square-free kernel $\delta=\delta(\varepsilon)$ of $N(1+\varepsilon)$ properly divides the discriminant $d(\varepsilon)$ of $\Q(\varepsilon)$ and takes on the value $+1$ at all the genus characters of $\Q(\varepsilon)$. This implies that $\sqrt{\varepsilon}=u\sqrt{\delta}+v\sqrt{d(\varepsilon)/\delta}$ for some rationals $u, v$.

\medskip
Now for (i),  if $\varepsilon=\varepsilon_{pq}$, then  clearly $\{\delta,d(\varepsilon)/\delta\}=\{p,q\}$ and thus in this case $\sqrt{\varepsilon_{pq}}=u\sqrt{p}+v\sqrt{q}$ for some $u,v\in \Q$.

\medskip
Next for (ii), if $\varepsilon=\varepsilon_{2p}$, then $d(\varepsilon)=8p$. Thus $\delta\in\{2,p,2p\}$, but $\delta\not=2p$ since otherwise $\Q(\varepsilon)=\Q(\sqrt{\varepsilon}\,)$. Thus $\sqrt{\varepsilon_{2p}}=u\sqrt{2}+v\sqrt{p}$, for some $u,v\in\Q$.

\medskip

Finally for (iii), if $\varepsilon=\varepsilon_p$, then $d(\varepsilon)=4p$. But then $\delta\in\{2,p,2p\}$. However, $\delta\not=p$, since one of the genus characters on $\Q(\sqrt{p}\,)$ is determined by the Kronecker symbol $(-4/\cdot)$ and $(-4/p)=-1$. Therefore, $\sqrt{\varepsilon_{p}}=u\sqrt{2p}+v\sqrt{2}$, for some $u,v\in\Q$.

\medskip
Now we will say more about the values of $u$ and $v$. In all these cases, let $K=\Q(\sqrt{\varepsilon}\,).$ Since $\sqrt{\varepsilon}$ is a unit, any norm applied to it is also a unit, and any trace is integral.

If $\varepsilon=\varepsilon_{pq}$, then $2u\sqrt{p}=\Tr_{K/\Q(\sqrt{p})}(\sqrt{\varepsilon_{pq}})\in \OO_{\Q(\sqrt{p})}=\Z[\sqrt{p}\,].$ Thus $a=2u\in \Z$.  Similarly, $b=2v\in\Z$. Moreover, $u^2p-v^2q=N_{K/\Q(\sqrt{p})}(\sqrt{\varepsilon_{pq}})$, which is a unit and rational, and therefore $=\pm 1$. This implies that
$a^2p-b^2q=\pm 4$ and thus $a\equiv b\bmod 2.$ Finally, notice that $\frac1{4}(a^2p-b^2q)=(p/q)$.

If $\varepsilon=\varepsilon_{2p}$, then again we see that $2u=x\in\Z$ and $2v=y\in\Z$ and moreover that $u^22-v^2p=\pm 1$, or equivalently, that $x^22-y^2p=\pm 4$. Suppose $x$ is odd. Then $2-y^2p\equiv 4\bmod 8$. Thus $y^2p\equiv -2\bmod 8$, which is impossible. This implies that $x$ is even, i.e., $u\in\Z$. Clearly, then $y$ is even, so $v\in\Z$. Thus
$a^22-b^2p=\pm 1$, with $a=u$ and $b=v$.  It follows easily that $a^22-b^2p=(2/p).$

Finally, suppose $\varepsilon=\varepsilon_p$, and so $\sqrt{\varepsilon_p}=u\sqrt{2p}+v\sqrt{2},$
for some rational $u,v$. But we now claim that $2u$ and $2v$ are odd integers. The same argument as above shows they are integers, so we need only consider their parity. But also as above we see that
$u^22p-v^22=\pm 1.$
This in turn easily implies that $a=2u$ and $b=2v$ are odd integers. Thus
$(a^2p-b^2)/2=\pm 1.$ Since $a^2p-b^2\equiv p-1\bmod 8$, we see that
$$\frac1{2}(a^2p-b^2)\equiv \frac{p-1}{2}\bmod 4\equiv  \left\{
                                  \begin{array}{ll}
                                    +1 \bmod 4, & \hbox{if}\;\;p\equiv 3\bmod 8, \\
                                    -1 \bmod 4, & \hbox{if}\;\;p\equiv 7\bmod 8,
                                  \end{array}
                                \right. $$
and therefore,
$$\frac1{2}(a^2p-b^2)=\left\{
               \begin{array}{ll}
                 +1, & \hbox{if}\;\;-p\equiv 5\bmod 8, \\
                 -1, & \hbox{if}\;\; -p\equiv 1\bmod 8,
               \end{array}
             \right.=-\big(\frac{-p}{2}\big).$$
\end{proof}

\begin{lem}\label{L7} Let $k$ be a quadratic number field satisfying the conditions in case $(c_2)$. Then

\medskip
$$\sqrt{\varepsilon_{k}}=\frac1{2}(a\sqrt{p_2}+b\sqrt{p_1p_3p_4})$$ for some  $a,b\in\Z$ with $a\equiv b\bmod 2.$
Moreover, $$\frac1{4}(a^2p_2-b^2p_1p_3p_4)=-1.$$
\end{lem}

\begin{proof} Since $\varepsilon_k$  has positive norm, we have (as in the proof of Lemma~\ref{L6}) $\sqrt{\varepsilon_k}=u\sqrt{\delta}+v\sqrt{d_k/\delta}$ for some rationals $u, v$. Furthermore, $\chi_\ell(\delta)=+1$ for all the genus characters $\chi_\ell$ of the field $k$. These characters can be obtained from the following table of Kronecker symbols given in Appendix~I. Let $\chi_i$ be the genus character associated with $d_i$, then

\medskip
\begin{tabular}{c|cccc}
   $\chi_i(p_j)$ & $\chi_1$ & $\chi_2$ & $\chi_3$ & $\chi_4$ \\ \hline
   \rsp $p_1$ & $-1 $ & $+1$ & $-1$ & $+1$  \\
   \rsp $p_2$ &  $-1 $ & $-1$ & $-1$ & $-1$  \\
   \rsp $p_3$ &  $+1 $ & $+1$ & $-1$ & $-1$  \\
   \rsp $p_4$ &  $-1 $ & $+1$ & $+1$ & $-1$  \\
     \end{tabular}.
\medskip

Therefore $\chi_i(p_1p_3p_4)=+1$ for all the genus characters, we have $\delta =p_1p_3p_4$, and thus
$\sqrt{\varepsilon_k}=u\sqrt{p_2}+v\sqrt{p_1p_3p_4}$ for some rational $u$ and $v$.
Notice that as in the proof of  Lemma~\ref{L6} we have $2u=a$, $2v=b$ for some $a,b\in\Z$ with $a\equiv b\bmod 2$ and therefore
$$\sqrt{\varepsilon_{k}}=\frac1{2}(a\sqrt{p_2}+b\sqrt{p_1p_3p_4}),$$
and moreover
$$\frac1{4}(a^2p_2-b^2p_1p_3p_4)=-(d_2/p_i)=-1,$$
 for $i=1,3,4$.
\end{proof}

By the lemmas, notice that $\sqrt{\varepsilon_{ij}\varepsilon_{\ell m}},\,\sqrt{\varepsilon_{ij}\varepsilon_k}\in E_{k^1}.$
Thus we can see that in particular
$$u_1=\sqrt{\varepsilon_{12}\varepsilon_{13}}, \; u_2=\sqrt{\varepsilon_{12}\varepsilon_{14}},\;u_3=\sqrt{\varepsilon_{12}\varepsilon_{23}},$$ $$ u_4=\sqrt{\varepsilon_{12}\varepsilon_{24}}, \; u_5=\sqrt{\varepsilon_{12}\varepsilon_{34}}, \; u_6=\sqrt{\varepsilon_{12}\varepsilon_k},$$
is a list of six independent units such that $(E_{k^1}:\la -1,u_1,u_2,u_3,u_4,u_5,u_6,\varepsilon_k\ra)=2$. Hence $E_{k^1}$  contains another unit $u_0$ and thus
$$E_{k^1}=\la -1,u_1,u_2,u_3,u_4,u_5,u_6,\varepsilon_k, u_0\ra.$$
We will determine the signs of $u_i^\sigma$ for $i=1,\dots,6$ and $\sigma\in\Gal(k^1/\Q)$, by looking in the field
$\mathcal{K}=\Q(\sqrt{p_1},\sqrt{p_2},\sqrt{p_3},\sqrt{p_4})\supseteq k^1$. We have $$\Gal(\mathcal{K}/\Q)=\la \tau_1,\tau_2, \tau_3, \tau_4\ra, \;\text{where}\;\; \tau_i(\sqrt{p_j}=\left\{
                                                              \begin{array}{ll}
                                                                -\sqrt{p_j}, &\;\; \hbox{if}\; j=i, \\
                                                                \phantom{-}\sqrt{p_j}, & \;\;\hbox{if not}.
                                                              \end{array}
                                                            \right.$$
Hence $\sigma_i=\tau_i|_{k^1}$ for $i=1,2,3$. We thus have the following table of some conjugates of some units in $\mathcal{K}$.

\medskip
\noindent\begin{tabular}{c|ccccccc}
 $\eta^\tau$ & $\tau_1$ & $\tau_2$ & $\tau_3$ & $\tau_1\tau_2$ & $\tau_1\tau_3$ & $\tau_2\tau_3$ & $\tau_1\tau_2\tau_3 $\\ \hline
\rsp $\sqrt{\varepsilon_{12}}$ & $-1/\sqrt{\varepsilon_{12}}$ & $1/\sqrt{\varepsilon_{12}}$ & $\sqrt{\varepsilon_{12}}$ & $-\sqrt{\varepsilon_{12}}$ & $-1/\sqrt{\varepsilon_{12}}$ & $1/\sqrt{\varepsilon_{12}}$ & $-\sqrt{\varepsilon_{12}}$\\
\rsp $\sqrt{\varepsilon_{13}}$ & $1/\sqrt{\varepsilon_{13}}$ & $\sqrt{\varepsilon_{13}}$ & $-1/\sqrt{\varepsilon_{13}}$ & $1/\sqrt{\varepsilon_{13}}$ & $-\sqrt{\varepsilon_{13}}$ & $-1/\sqrt{\varepsilon_{13}}$ & $-\sqrt{\varepsilon_{13}}$\\
\rsp $\sqrt{\varepsilon_{14}}$ & $-1/\sqrt{\varepsilon_{14}}$ & $\sqrt{\varepsilon_{14}}$ & $\sqrt{\varepsilon_{14}}$ & $-1/\sqrt{\varepsilon_{14}}$ & $-1/\sqrt{\varepsilon_{14}}$ & $\sqrt{\varepsilon_{14}}$ & $-1/\sqrt{\varepsilon_{14}}$\\
\rsp $\sqrt{\varepsilon_{23}}$ & $\sqrt{\varepsilon_{23}}$ & $1/\sqrt{\varepsilon_{23}}$ & $-1/\sqrt{\varepsilon_{23}}$ & $1/\sqrt{\varepsilon_{23}}$ & $-1/\sqrt{\varepsilon_{23}}$ & $-\sqrt{\varepsilon_{23}}$ & $-\sqrt{\varepsilon_{23}}$\\
\rsp $\sqrt{\varepsilon_{24}}$ & $\sqrt{\varepsilon_{24}}$ & $1/\sqrt{\varepsilon_{24}}$ & $\sqrt{\varepsilon_{24}}$ & $1/\sqrt{\varepsilon_{24}}$ & $\sqrt{\varepsilon_{24}}$ & $1/\sqrt{\varepsilon_{24}}$ & $1/\sqrt{\varepsilon_{24}}$\\
\rsp $\sqrt{\varepsilon_{34}}$ & $\sqrt{\varepsilon_{34}}$ & $\sqrt{\varepsilon_{34}}$ & $1/\sqrt{\varepsilon_{34}}$ & $\sqrt{\varepsilon_{34}}$ & $1/\sqrt{\varepsilon_{34}}$ & $1/\sqrt{\varepsilon_{34}}$ & $1/\sqrt{\varepsilon_{34}}$\\
\rsp $\sqrt{\varepsilon_{k}}$ & $-1/\sqrt{\varepsilon_{k}}$ & $1/\sqrt{\varepsilon_{k}}$ & $-1/\sqrt{\varepsilon_{k}}$ & $-\sqrt{\varepsilon_{k}}$ & $\sqrt{\varepsilon_{k}}$ & $-\sqrt{\varepsilon_{k}}$ & $1/\sqrt{\varepsilon_{k}}$\\
\end{tabular}
\medskip

For example,
 $$\sqrt{\varepsilon_{12}}^{\,\tau_1}=(a_{12}\sqrt{p_2}+b_{12}\sqrt{p_1})^{\tau_1}=
a_{12}\sqrt{p_2}-b_{12}\sqrt{p_1}=-(d_2/p_1)/\sqrt{\varepsilon_{12}}=-1/\sqrt{\varepsilon_{12}}\,,$$
by the proof of Lemma~\ref{L6} and the value of the genus characters above.

From this table we get the following signs of the conjugates of the units $u_i$.

\medskip
\begin{tabular}{c|ccccccc}
 sign\,$\eta^\sigma $ & $\sigma_1$ & $\sigma_2$ & $\sigma_3$ & $\sigma_1\sigma_2$ & $\sigma_1\sigma_3$ & $\sigma_2\sigma_3$ & $\sigma_1\sigma_2\sigma_3 $\\ \hline
\rsp $u_1$ & $-1$ & $+1$ & $-1$ & $-1$  & $+1$ & $-1$ & $+1$ \\
\rsp $u_2$ & $+1$ & $+1$ & $+1$ & $+1$  & $+1$ & $+1$ & $+1$ \\
\rsp $u_3$ & $-1$ & $+1$ & $-1$ & $-1$  & $+1$ & $-1$ & $+1$  \\
\rsp $u_4$ & $-1$ & $+1$ & $+1$ & $-1$  & $-1$ & $+1$ & $-1$ \\
\rsp $u_5$ &  $-1$ & $+1$ & $+1$ & $-1$  & $-1$ & $+1$ & $-1$ \\
\rsp $u_6$ & $+1$ & $+1$ & $-1$ & $+1$  & $-1$ & $-1$ & $-1$ \\
\rsp $\varepsilon_k $ & $+1$ & $+1$ & $+1$ & $+1$  & $+1$ & $+1$ & $+1$ \\ \hline
\rsp $S$ & $ \uparrow$      &      & $ \uparrow$   & $ \uparrow$   &  & $\uparrow $ & \\
\end{tabular}

\medskip
Here, for ease of reference, we have added a row showing by arrows the relevant columns for each of the sets $S$. From the table we thus see that the unit group $H$ satisfies:
$$H\supseteq\la -u_1,u_2,u_1u_3,u_4u_5,u_3u_4u_6,\varepsilon_k, E_{k^1}^2\ra,$$
and we observe that  $(\la -u_1,u_2,u_1u_3,u_4u_5,u_3u_4u_6,\varepsilon_k, E_{k^1}^2\ra:E_{k^1}^2)=2^6$. Since \\ $(E_{k^1}:E_{k^1}^2)=2^8,$ we have $(E_{k^1}:H)\leq 4.$ Thus $\#\Am_2(K/k^1)\geq 2$ and therefore $h_2(K)\geq 2$, as desired for $(c_2)$.

\bigskip

Now we consider $(d_5),(d_8)$.  In these cases,
$d_k=-4d_1d_2d_3=4p_1p_2p_3$ and
$$ k^1=\Q(\sqrt{p_1},\sqrt{p_2},\sqrt{p_3}\,),\qquad k_+^1=k^1(i).$$
Let $F=\Q(\sqrt{d_1d_2d_3}\,)=\Q(\sqrt{-p_1p_2p_3}\,)$; then $d_F=d_1d_2\cdot d_3$ is a $C_4$-splitting of $d_F$. This implies that there is an unramified cyclic quartic extension $\cF$ of $F$. Thus there is a primitive rational solution $(a,b,c)$ satisfying $a^2=b^2d_3+c^2d_1d_2$ such that $\alpha=a+c\sqrt{d_1d_2}>0$ with $\cF=F(\sqrt{\alpha}\,)$. Moreover, $\alpha\alpha'=a^2-c^2d_1d_2=b^2d_3$.

Now consider the extension $\cF/\Q(\sqrt{d_1d_2}\,).$ Then by shifting by $k^1$ we obtain the finitely unramified $V_4$-extension $L/k^1$ with $L=k^1\cF=k^1(\sqrt{-1},\sqrt{\alpha}\,)$:
$$L=k^1(i,\sqrt{\alpha}\,)$$
$$k_+^1=k^1(i),\qquad K=k^1(\sqrt{\alpha}\,),\quad K'=k^1(\sqrt{\alpha'}\,)=k^1(\sqrt{d_3\alpha}\,)$$
$$k^1$$
As above it suffices to show that $h_2(K)\geq 2.$  To this end we apply the Ambiguous Class Number Formula to the extension $K/k^1$,
$$2^{\rank\Cl_2(K)}=\#\Am_2(K/k^1)=\frac{2^{t-1}}{(E_{k^1}:H_{K/k^1})},$$
where $t=t_\infty$ is the number of infinite places of $k^1$  that ramify in $K$. We can see that  $t=t_\infty=\#\{\sigma\in\Gal(k^1/\Q):\alpha^\sigma<0\}$. Now $\Gal(k^1/\Q)=\la \sigma_1,\sigma_2,\sigma_3\ra$ where
$$\sigma_i(\sqrt{p_j}\,)=\left\{
                              \begin{array}{ll}
                                -\sqrt{p_j}, & \; \hbox{if}\;\; j=i, \\
                                \phantom{-}\sqrt{p_j}, & \; \hbox{if}\;\; j\not=i.
                              \end{array}
                            \right.$$
Since $\Q(\alpha)=\Q(\sqrt{p_1p_2}\,)$ we see that
$$\alpha^\sigma <0\;\;\text{iff}\;\;\sigma\in S=\{\sigma_1,\sigma_2, \sigma_1\sigma_3,\sigma_2\sigma_3\}.$$
Thus $t=4$ and so $2^{t-1}=2^3.$

Now consider $H.$ By local-global properties, we have
$$H=\{\varepsilon\in E_{k^1}:\varepsilon^\sigma>0, \;\sigma\in S\}.$$

We now  compute $E_{k^1}$ As before $q(k^1/\Q)=2^7$, where $q(k^1/\Q)=(E_{k^1}:E)$, with $E$ the  subgroup of $E_{k^1}$ generated by all the units in the quadratic subfields of $k^1$, i.e.,
$$E=\la -1,\varepsilon_{12},\varepsilon_{13},\varepsilon_{14},\varepsilon_{23},\varepsilon_{24},\varepsilon_{34},\varepsilon_k \ra=\la -1,\varepsilon_{12},\varepsilon_{13},\varepsilon_{23},\varepsilon_1,\varepsilon_{2}, \varepsilon_{3},\varepsilon_k \ra$$
 Notice that $\varepsilon_{j4}$ is the fundamental unit in $\Q(\sqrt{p_j}\,),$ thus to simply notation a bit, we let $\varepsilon_{j4}=\varepsilon_j$.

We also have the following result analogous to Lemma~\ref{L7}.
\begin{lem}\label{L8} Let $k$ be a quadratic number field satisfying the conditions in cases $(d_5)$ or $(d_8)$. Then
$$\sqrt{\varepsilon_{k}}=\frac1{2}(a\sqrt{p_3}+b\sqrt{p_1p_2})$$ for some  $a,b\in\Z$ such that $a\equiv b\bmod 2.$
Moreover, $$\frac1{4}(a^2p_3-b^2p_1p_2)=-1.$$
\end{lem}

\begin{proof} By the  Kronecker symbols given in Appendix~I and where  $\chi_i$ is the genus character associated with $d_i$, then

\medskip
\begin{tabular}{c|cccc}
   $\chi_i(p_j)$ & $\chi_1$ & $\chi_2$ & $\chi_3$ & $\chi_4$ \\ \hline
   \rsp $p_1$ & $-1 $ & $+1$ & $+1$ & $-1$  \\
   \rsp $p_2$ &  $-1 $ & $+1$ & $+1$ & $-1$  \\
   \rsp $p_3$ &  $-1 $ & $-1$ & $-1$ & $-1$  \\
   \rsp $p_4$ &  $\pm 1 $ & $-1$ & $+1$ & $\mp 1$  \\
     \end{tabular}.

\noindent Here for $\pm 1$ and $\mp 1$, the top sign occurs in the case $(d_5)$ and the bottom for $(d_8)$.
\medskip

Therefore $\chi_i(p_1p_2)=+1$ for all the genus characters, we have $\delta =p_1p_2$, and thus
$\sqrt{\varepsilon_k}=u\sqrt{p_3}+v\sqrt{p_1p_2}$ for some rational $u$ and $v$.

The rest of the proof again follows similarly to that in Lemma~\ref{L6}.   Moreover, notice that
$$\frac1{4}(a^2p_3-b^2p_1p_2)=-(d_3/p_i)=-1$$
 for $i=1,2$.
\end{proof}

By Lemmas~\ref{L6} and \ref{L8}, we can see that $\sqrt{\varepsilon_{12}},\sqrt{\varepsilon_{13}},\sqrt{\varepsilon_{23}},\sqrt{\varepsilon_{k}}\in k^1,$ as well as
$\sqrt{\varepsilon_1\varepsilon_2},\sqrt{\varepsilon_1\varepsilon_3},\sqrt{\varepsilon_2\varepsilon_3}.$
Thus we  see that in particular
$$u_1=\sqrt{\varepsilon_{12}}, \; u_2=\sqrt{\varepsilon_{13}},\;u_3=\sqrt{\varepsilon_{23}},$$ $$ u_4=\sqrt{\varepsilon_{k}}, \; u_5=\sqrt{\varepsilon_{1}\varepsilon_{2}}, \; u_6=\sqrt{\varepsilon_{1}\varepsilon_3},$$
is a list of six independent units such that $(E_{k^1}:\la -1,u_1,u_2,u_3,u_4,u_5,u_6,\varepsilon_1\ra)=2$. Hence $E_{k^1}$  contains another unit $u_0$ and thus
$$E_{k^1}:\la -1,u_1,u_2,u_3,u_4,u_5,u_6,\varepsilon_1, u_0\ra.$$
We can  determine the signs of $u_i^\sigma$ for $i=1,\dots,6$ similarly to above and see that

\medskip
\begin{tabular}{c|ccccccc}
 sign\,$\eta^\sigma $ & $\sigma_1$ & $\sigma_2$ & $\sigma_3$ & $\sigma_1\sigma_2$ & $\sigma_1\sigma_3$ & $\sigma_2\sigma_3$ & $\sigma_1\sigma_2\sigma_3 $\\ \hline
\rsp $u_1$ & $-1$ & $+1$ & $+1$ & $-1$  & $-1$ & $+1$ & $-1$ \\
\rsp $u_2$ & $-1$ & $+1$ & $+1$ & $-1$  & $-1$ & $+1$ & $-1$\\
\rsp $u_3$ & $+1$ & $-1$ & $+1$ & $-1$  & $+1$ & $-1$ & $-1$  \\
\rsp $u_4$ & $-1$ & $-1$ & $+1$ & $+1$  & $-1$ & $-1$ & $+1$ \\
\rsp $u_5$ &  $\pm 1$ & $-1$ & $+1$ & $\mp 1$  & $\pm 1$ & $-1$ & $\mp 1$ \\
\rsp $u_6$ & $\pm 1$ & $+1$ & $+1$ & $\pm 1$  & $\pm 1$ & $+1$ & $\pm 1$ \\
\rsp $\varepsilon_1 $ & $+1$ & $+1$ & $+1$ & $+1$  & $+1$ & $+1$ & $+1$ \\ \hline
\rsp $S$ & $ \uparrow$      & $ \uparrow$     &    &   & $ \uparrow$   & $\uparrow $ & \\
\end{tabular}

\medskip
Here, for ease of reference, we have added a row showing by arrows the relevant columns for each of the sets $S$. From the table we thus see that the unit group $H$ satisfies:
$$H\supseteq \left\{
               \begin{array}{ll}
                 \la u_1u_2,u_3u_5u_6,-u_4,u_6, -u_2u_3,\varepsilon_1, E_{k^1}^2\ra, & \hbox{for}\;\;(d_5), \\
                  \la u_1u_2,u_3u_5u_6,-u_4,-u_5, -u_2u_3,\varepsilon_1, E_{k^1}^2\ra, & \hbox{for}\;\;(d_8).
               \end{array}
             \right.
$$
and so  $(E_{k^1}:H)\leq 4.$ Thus $\#\Am_2(K/k^1)\geq 2$ and therefore $h_2(K)\geq 2$, as desired.

\bigskip

\noindent{\em Examples.} We now give examples of fields $k$ in each of the cases $(c_2), (d_5),(d_8)$. Let $K=k^1(\sqrt{\alpha}\,)$ and we will see that $h_2(K)\geq 2$ and so $k$ has narrow $2$-class field tower of length at least $3$.

\medskip
For $(c_2)$ let
$$d_k=7\cdot 3\cdot 43\cdot 31,\quad \Cl_2(k_+^1)\simeq (2,8,8), \quad \alpha=383+26\sqrt{7\cdot 31},\quad h_2(K)=4.$$

\medskip
For $(d_5)$ let
$$d_k=7\cdot 3\cdot 47\cdot 4,\quad \Cl_2(k_+^1)\simeq (2,4,4), \quad \alpha=17+4\sqrt{21},\quad h_2(K)=2.$$

\medskip
For $(d_8)$ let
$$d_k=11\cdot 43\cdot 7\cdot 4,\quad \Cl_2(k_+^1)\simeq (2,4,4), \quad \alpha=13+2\sqrt{44},\quad h_2(K)=2.$$

\bigskip

All of this establishes the Main Theorem.

\bigskip
We finally consider our last cases.

\bigskip
\section{{\bf (D) the cases for which $G^+/G^+_3\simeq 64.150$}}

We have two cases to consider, $(c_3)$ and $(d_1)$. As easily seen in Appendix I, these involve precisely the fields $k=\Q(\sqrt{d_1d_2d_3d_4}\,)$, where $d_j$ are distinct negative prime discriminants each divisible by the prime $p_j$ such that

\medskip

(1) if $d_k\equiv 4\bmod 8$, then $d_4=-4$,

(2) $(d_1/p_2)=(d_2/p_3)=(d_3/p_1)=-1$, $(d_1/p_4)=(d_2/p_4)=(d_3/p_4)=+1.$

\medskip

Let's review what we know so far: Let $k^1$ and $k_+^1$ denote the Hilbert $2$-class field, resp., narrow $2$-class field, of $k$. Hence we have $k^1=\Q(\sqrt{d_1d_4},\sqrt{d_2d_4},\sqrt{d_3d_4}\,)$ with $\Cl_2(k)\simeq \Gal(k^1/k)\simeq (2,2)$, such that $h_2(k^1)=1$  and
$k_+^1=\Q(\sqrt{d_1},\sqrt{d_2},\sqrt{d_3}, \sqrt{d_4}\,)$ where $\Cl_2^+(k)\simeq \Gal(k_+^1/k)\simeq (2,2,2)$. All of this information can be found  in   \cite{BS3}. Moreover, we will see below that $\Cl_2(k_+^1)\simeq \Gal(k_+^2/k_+^1)$ has rank $3$.  Furthermore,
if we let $G^+=\Gal(k_+^2/k)$, then we have $G^+/G^+_3\simeq 64.150$. Our goal was to show  that the narrow 2-class field tower of $k$ terminates at $k_+^2$. However, we have not been able to prove this. Nevertheless, we will show below that if $\Cl_2(k_+^1)\subseteq (4,4,4)$, then the narrow $2$-class field of $k$ terminates at $k_+^2$.

We start with some group theoretical results and their applications. We have the following result about number fields, which is an immediate corollary  of the group theoretic result below.

\begin{prop}\label{P9} Let $k$ be a number field with $\rank\Cl_2(k)=3.$ If there are  distinct unramified quadratic extensions $K_1,K_2,K_3$ of $k$ with $K_3\subseteq K_1K_2$ and such that $h_2(K_j)=\frac1{2}h_2(k)$, $j=1,2,3$, then for all seven unramified quadratic extensions $K$ of $k$ we have $h_2(K)=\frac1{2}h_2(k)$. Moreover, the $2$-class field tower of $k$ terminates at $k^1$, i.e., $k^1=k^2$. Conversely, if $k^1=k^2$, then $h_2(K)=\frac1{2}h_2(k)$ for all unramified quadratic extensions of $k$.
\end{prop}

This proposition follows immediately from the group theoretic result:

\begin{prop}\label{P10} Let $G$ be a finite $2$-group such that $\rank G/G'=3$. Suppose $H_1,H_2,H_3$ are distinct maximal subgroups of $G$ with $H_1\cap H_2\subseteq H_3$ and such that $H_j'=G'$, for $j=1,2,3$. Then $G_2=G'=1$ and thus $H'=G'$ for all subgroups $H$ of $G$. Conversely, if $G'=1$, then (clearly!) $H'=G'$ for all subgroups $H$ of $G$.
\end{prop}
\begin{proof} Let $H_1,H_2,H_3$ be distinct maximal subgroups of $G$ with $H_1\cap H_2\subseteq H_3$ such that $H_j'=G'$ for $j=1,2,3.$ Then without loss of generality we choose generators $a_1,a_2,a_3$ of $G$ such that
$$H_1=\la a_1,a_2,G^2\ra, \;\; H_2=\la a_1,a_3,G^2\ra,\;\; H_3=\la a_1,a_2a_3,G^2\ra.$$
Hence, since we are assuming $G'=H_j'$\,,
$$G'=\la c_{12},c_{13},c_{23},G_3\ra,\; H_1'=\la c_{12},G_2^2G_3\ra,\; H_2'=\la c_{13},G_2^2G_3\ra,\; H_3'=\la c_{12}c_{13},G_2^2G_3\ra,$$
where $c_{ij}=[a_i,a_j]$.
We want to show that $G'=1$. To this end, we may assume that $G'^2G_3=1$. To see this, notice that
$$G'=1\Leftrightarrow  G'=G_3 \Leftrightarrow(G/G_3)'=1\Leftrightarrow (G/G_3)'=(G/G_3)'^{\,2} \Leftrightarrow $$
$$G'/G_3=G'^2G_3/G_3 \Leftrightarrow G'=G'^2G_3\Leftrightarrow (G/G'^2G_3)'=1.$$

With this assumption we have $G'=\la c_{12}, c_{13},c_{23}\ra$ with $c_{ij}^2=1$.  Moreover, $G'=H_1'=\la c_{12}\ra=H_2'=\la c_{13}\ra=H_3'=\la c_{12}c_{13}\ra.$ Therefore, $c_{12}=c_{13}=c_{12}c_{13}$ and thus in particular  $c_{12}=1$, whence  $G'=1.$
\end{proof}

Next, we will prove some results about $2$-groups $G$ for which $G/G_3$ is of type $64.150$.  Recall from \cite{HS} that if $G$ is a finite $2$-group with $G/G_3\simeq 64.150$, then
$$G=\la a_1,a_2,a_3: a_1^2\equiv c_{12},a_2^2\equiv c_{23},a_3^2\equiv c_{13}\bmod G_3\ra.$$

As promised above, we have the following proposition.
\begin{prop}\label{P11} Let $G$ be a finite $2$-group for which $G/G_3\simeq 64.150$ as presented above. Then
$$G'=\la a_1^2,a_2^2,a_3^2\ra=\la c_{12},c_{23},c_{13}\ra$$
and hence of rank $3$. Moreover, for the $j$-th term $G_j$ in the lower central series, we have
$$ G_j=\la a_1^{2^{j-1}}, a_2^{2^{j-1}}, a_3^{2^{j-1}}\ra=G_2^{2^{j-2}}.$$
\end{prop}
\begin{proof} We have
$$G'=\la a_1^2,a_2^2,a_3^2,G_3\ra=\la c_{12},c_{23},c_{13},G_3\ra,$$
and thus
$$G_3=\la c_{12i},c_{23i},c_{13i}, G_4: i=1,2,3\ra,$$
where $c_{ij\ell}=[c_{ij},a_\ell]$.
Then it follows that
$$c_{121}\equiv 1,\; c_{122}\equiv c_{12}^2,\; c_{123}\equiv 1\bmod G_4,$$
$$c_{231}\equiv 1,\; c_{232}\equiv 1,\; c_{233}\equiv c_{23}^2\bmod G_4,$$
$$c_{131}\equiv c_{13}^2,\; c_{132}\equiv 1,\; c_{133}\equiv 1\bmod G_4.$$
For example,
$$1=[a_1^2,a_1]\equiv [c_{12},a_1]=c_{121}\bmod G_4$$
and then
$$c_{12}^2c_{121}\equiv [a_1^2,a_2]\equiv [c_{12},a_2]=c_{122}\bmod G_4,$$
thus $c_{122}\equiv c_{12}^2\bmod G_4$. (Here and in various other places,  we are using several applications of Hilfssatz 1.2 in Chapter III of \cite{Hup}.)

Therefore, $G_3=\la c_{12}^2,c_{23}^2,c_{13}^2,G_4\ra=\la a_1^4,a_2^4,a_3^4, G_4\ra$ (since $a_{12}^2\equiv c_{12}\bmod G_3$ implies that $a_1^4\equiv c_{12}^2\bmod G_4$, etc.), whence $G'=\la a_1^2,a_2^2,a_3^2,G_4\ra$.
Next, we have
$$G_4=\la [c_{12}^2,a_i],[c_{23}^2,a_i],[c_{23}^2,a_i],G_5: i=1,2,3\ra.$$
Now notice that $[c_{12}^2,a_i]\equiv c_{12i}^2[c_{12i},c_{12}]\equiv c_{12i}^2\bmod G_5,$ since $[c_{123},c_{12}]\in G_5$, for recall that $[G_i,G_j]\subseteq G_{i+j}$, cf.,\,\cite{Hup}. But $c_{12i}^2\equiv (1,c_{12}^4,1)\bmod G_5$ for $i=(1,2,1)$, etc.  Therefore $G_4=\la c_{12}^4,c_{23}^4,c_{13}^4,G_5\ra=\la a_1^8, a_2^8, a_3^8, G_5\ra,$ and so
$G'=\la a_1^2,a_2^2,a_3^2,G_5\ra =\la c_{12},c_{23},c_{13}, G_5\ra.$ Continuing  we get
$G'=\la a_1^2,a_2^2,a_3^2,G_j\ra =\la c_{12},c_{23},c_{13}, G_j\ra$, for all $j$. Since $G_j=1$ for sufficiently large $j$,
we have $G'=\la a_1^2,a_2^2,a_3^2\ra =\la c_{12},c_{23},c_{13}\ra$. Moreover, $G_j=G_2^{2^{j-2}}$, as desired.
\end{proof}

For $G/G_3\simeq 64.150$ the following proposition gives a sufficient condition  for $G$ to be metabelian.
\begin{prop}\label{P12} Let $G$ be a finite $2$-group for which $G/G_3\simeq 64.150$ and such that $8$-$\rank G'/G''=0$, i.e., $G'/G''\subseteq (4,4,4).$ Then $G$ is metabelian, i.e., $G''=1$.
\end{prop}
\begin{proof} By the proof of Proposition~\ref{P10} applied to $G'$, in order to show $G''=1$, we may assume $(G'')^2(G')_3=1.$ Assuming this and since $G'=\la c_{12},c_{13},c_{23}\ra$, we have
$$G''=\la [c_{12},c_{13}], [c_{12},c_{23}], [c_{13},c_{23}]\ra,$$
and for which each $[c_{ij},c_{\ell m}]^2=1$.

By our hypothesis $(G')^4\subseteq G''$. By the previous proposition, $(G')^4=G_4$ and thus $G_4\subseteq G''$. But for any $p$-group, $G''\subseteq G_4$, see \cite{Hup}, and thus $G''=G_4.$ With this, and Proposition~\ref{P11}, we  have $G_5=G_4^2=(G'')^2=1$. We will now show that  $[c_{ij},c_{\ell m}]\equiv 1\bmod G_5$ and thus $=1$. To see this, we have for instance
$$[c_{12},c_{13}]=[a_1^2,a_3^2]=[a_1^2,a_3]^2[a_1^2,a_3,a_3]=(c_{13}^2c_{131})^2c_{133}^2c_{1313}\equiv c_{1313}\bmod G_5,$$
since $c_{133}, c_{13}^2c_{131}\equiv 1\bmod G_4$ (and so their squares are in $G_5$). Moreover,
$$c_{1313}=[c_{131},a_3]\equiv [c_{13}^2,a_3]\equiv c_{133}^2[c_{13},a_3,c_{13}]\equiv 1\bmod G_5.$$
Therefore, $[c_{12},c_{13}]=1$. The other generators of $G''$ are also $=1$ in a similar fashion.

Therefore, $G''=1$, as desired.
\end{proof}

We thus immediately get the following corollary.

\begin{cor}\label{C1} Let $k$ be a real quadratic number field with discriminant $d_k=d_1d_2d_3d_4$, such that the $d_j$ are negative prime discriminants divisible by the primes $p_j$  (respectively) satisfying

(1) if $d_k\equiv 4\bmod 8$, then $d_4=-4$,

(2) $(d_1/p_2)=(d_2/p_3)=(d_3/p_1)=-1$, $(d_1/p_4)=(d_2/p_4)=(d_3/p_4)=+1.$

\medskip

\noindent If  $8$-$\rank\Cl_2(k_+^1)=0$, then the narrow $2$-class field tower of $k$ terminates at $k_+^2$.

\end{cor}

\bigskip
We now show a method for determining whether or not the narrow $2$-class field tower of $k$ terminates at $k_+^2$ for any $k$ for which $G^+/G^+_3\simeq 64.150.$

\subsection{Construction of some quartic extensions}
Let $$F_1=\Q(\sqrt{d_1d_3d_4}\,),\;F_2=\Q(\sqrt{d_1d_2d_4}\,),\;F_3=\Q(\sqrt{d_2d_3d_4}\,).$$ Their discriminants have $C_4$-splittings: $$d_{F_1}=d_1\cdot d_3d_4,\;d_{F_2}=d_2\cdot d_1d_4,\;d_{F_3}=d_3\cdot d_2d_4.$$
We will only need to examine two of these fields, say $F_1$ and $F_2$. From the $C_4$-splittings we obtain unramified cyclic quartic extensions $\mathcal{F}_j$ of $F_j$ ($j=1,2$). If we let $L_j=k^1\mathcal{F}_j$, then we have unramified quadratic extensions of $k_+^1$  such that the composite $L_1L_2$ is a $V_4$-extension of $k_+^1$ and hence containing a third quadratic extension $L_0$ of $k_+^1$. Since the rank of $\Cl_2(k_+^1)$ is $3$, if we can show that $h_2(L_j)=\frac1{2}h_2(k_+^1)$, then by Proposition~\ref{P9}, we will have our desired result that $k_+^2=k_+^3$.

Now back to constructing $\mathcal{F}_j$. By the $C_4$-splittings of the  $d_{F_j}$, there are primitive rational solutions $(a_j,b_j,c_j)$ satisfying
$$a_1^2=d_1b_1^2+d_3d_4c_1^2,\quad  a_2^2=d_2b_2^2+d_1d_4c_2^2$$
such that if we let
$$\alpha_1=a_1+c_1\sqrt{d_3d_4},\quad \alpha_2=a_2+c_2\sqrt{d_1d_4},$$
then $\mathcal{F}_j=F_j(\sqrt{\alpha_j}\,)$ are our desired extensions, cf.\,\cite{Lem}.  We'll now consider $\mathcal{F}_1$ in detail.
Notice that $\Q(\alpha_1)=\Q(\alpha_1')=\Q(\sqrt{d_3d_4}\,)$ and that \\ $\alpha_1\alpha_1'=a_1^2-c_1^2d_3d_4=b_1^2d_1<0$. Hence without loss of generality, assume $\alpha_1>0$ and $\alpha_1'<0$. Now, we have
$\sqrt{\alpha_1'}=\pm b_1\sqrt{d_1\alpha_1}/{\alpha_1},$
thus $\Q(\sqrt{d_3d_4},\sqrt{\alpha_1'}\,)=\Q(\sqrt{d_3d_4},\sqrt{d_1\alpha_1}\,).$ Finally observe that
$\mathcal{F}_1=F_1(\sqrt{\alpha_1}\,)=\Q(\sqrt{d_3d_4},\sqrt{d_1}, \sqrt{\alpha_1}\,).$

\subsection{Three useful quadratic extensions of $k_+^1$}

Shifting the extension \\$\mathcal{F}_1/\Q(\sqrt{d_3d_4}\,)$ by $k^1$ yields the finitely unramified $V_4$-extension $L_1=k^1(\sqrt{d_1},\sqrt{\alpha_1}\,)$ of $k^1$ for which we have $L_1=k_+^1(\sqrt{\alpha_1}\,)$ is an (everywhere) unramified quadratic extension of $k_+^1$.

The same type of argument shows that $L_2=k^1(\sqrt{d_2},\sqrt{\alpha_2}\,)$ is also a finitely unramified $V_4$-extension of $k^1$ for which $L_2=k_+^1(\sqrt{\alpha_2}\,)$ is an unramified quadratic extension of $k_+^1$.  Moreover, the composite $L_1L_2$ is an (everywhere) unramified $V_4$-extension of $k_+^1$ with $L_0=k_+^1(\sqrt{\alpha_1\alpha_2}\,)$, the third quadratic extension of $k_+^1$ in $L_1L_2$. By Proposition~\ref{P9}, we have the following result:

\begin{prop}\label{P13} The narrow $2$-class field tower of $k$ terminates at $k_+^2$ if and only if
$h_2(L_j)/h_2(k_+^1)=1/2$, for $j=0,1,2$.
\end{prop}

\subsection{Two examples} We now consider two examples of fields $k$, one in case $(c_3)$ and the other in $(d_1)$. We will use the above information to show that the narrow $2$-class field tower for each terminates at $k_+^2$.

\medskip
\noindent {\em Example 1.}  Let $k$ satisfy $d_k=3\cdot 8\cdot 11\cdot 23$. Then $k$ is in case $(c_3)$ and thus $G^+/G^+_3\simeq 64.150$.  By Pari, we obtained $\Cl_2(k_+^1)\simeq (2,4,4)$, and therefore by Corollary~\ref{C1} we conclude that the narrow $2$-class field tower of $k$ has length $2$.

We will now corroborate this using Proposition~\ref{P13}. Let $L_j=k_+^1(\sqrt{\alpha_j}\,)$ for $j=0,1,2$, where
$$ \alpha_1=-1+\sqrt{24},\quad \alpha_2=5+\sqrt{33},\quad \alpha_0=\alpha_1\alpha_2.$$
Then by Pari, we obtained
$$\Cl_2(L_1)\simeq (4,4),\quad \Cl_2(L_2)\simeq (4,4),\quad \Cl_2(L_0)\simeq (2,2,4).$$
Thus we see  by Proposition~\ref{P13} that $k$ has narrow $2$-class field tower of length $2$.

\medskip
\noindent {\em Example 2.}  Let $k$ satisfy $d_k=7\cdot 31\cdot 23\cdot 4$. Then $k$ is in case $(d_1)$ and thus again $G^+/G^+_3\simeq 64.150$.  By Pari, we obtained $\Cl_2(k_+^1)\simeq (2,4,8)$, and therefore we cannot use  Corollary~\ref{C1}.  However, we can use Proposition~\ref{P13}. Let $L_j=k_+^1(\sqrt{\alpha_j}\,)$ for $j=0,1,2$, where
$$ \alpha_1=19+2\sqrt{92},\quad \alpha_2=9+2\sqrt{28},\quad \alpha_0=\alpha_1\alpha_2.$$
Then by Pari, we obtained
$$\Cl_2(L_1)\simeq (2,4,4),\quad \Cl_2(L_2)\simeq (4,8),\quad \Cl_2(L_0)\simeq (4,8).$$
Thus we see that $k$ has narrow $2$-class field tower of length $2$.

\newpage

\begin{center} {\bf Appendix~I} \end{center}

\bigskip

See the description as the end of the appendix, also see \cite{BS3}:

\noindent{\bf Type I}\;\;  ($d_1,d_2>0$, $d_3,d_4<0$, no $d_j=-4$)

For each case, $(d_3/p_4)=+1$ and

\medskip

\begin{tabular}{c|ccccc|c|c}
   \text{case} & $(d_1/p_2)$ & $(d_1/p_3)$ & $(d_2/p_3)$
            & $(d_1/p_4)$ & $(d_2/p_4)$ & $G$ & $G^+/G_3^+$ \\ \hline
   \rsp $(a_1)$ & $+1$ & $+1$ & $+1$ & $-1$ & $-1$ & $Q_g$\;or\, $D$  &  64.144 \\
   \rsp $(a_2)$ & $+1$ & $-1$ & $-1$ & $+1$ & $+1$ & $Q_g$\;or\, $D$  &  64.144\\
   \rsp $(a_3)$ & $-1$ & $+1$ & $-1$ & $+1$ & $+1$  & $Q_g$\;or\, $D$  &  64.144  \\
   \rsp $(a_4)$ & $-1$ & $+1$ & $+1$ & $+1$ & $-1$ & $Q_g$\;or\, $D$  &  64.144 \\
   \rsp $(a_5)$ & $+1$ & $-1$ & $+1$ & $-1$ & $-1$  & $D$ & 32.033\\
   \rsp $(a_6)$ & $+1$ & $-1$ & $-1$ & $+1$ & $-1$  & $D$ & 32.033\\
   \rsp $(a_7)$ & $+1$ & $-1$ & $+1$ & $+1$ & $-1$ & $D$ & 64.144\\
   \rsp $(a_8)$ & $-1$ & $+1$ & $-1$ & $+1$ & $-1$  & $S$\;or\, $D$  &  32.033  \\
   \rsp $(a_9)$ & $-1$ & $+1$ & $+1$ & $-1$ & $-1$ & $Q$ & 64.147 \\
   \rsp $(a_{10})$ & $-1$ & $-1$ & $-1$ & $+1$ & $+1$ & $Q$ & 64.147\\
   \rsp $(a_{11})$ & $-1$ & $-1$ & $+1$ & $-1$ & $-1$ & $V_4$  &  32.039 \\
   \rsp $(a_{12})$ & $-1$ & $-1$ & $-1$ & $+1$ & $-1$ & $V_4$  & 32.039 \\
   \rsp $(a_{13})$ & $-1$ & $-1$ & $+1$ & $+1$ & $-1$  & $V_4$  & 32.034 \\
  \end{tabular}

\bigskip

\noindent{\bf Type II}\;\;  ($d_1,d_2>0$, $d_3,d_4<0$,  $d_4=-4$)

\medskip
For each case, $(d_3/2)=\pm 1$ and

\medskip

\begin{tabular}{c|ccccc|c|c}
   \text{case} & $(d_1/p_2)$ & $(d_1/p_3)$ & $(d_2/p_3)$
            & $(d_1/2)$ & $(d_2/2)$ & $G$ & $G^+/G_3^+$ \\ \hline
   \rsp $(b_1)$ & $+1$ & $-1$ & $-1$ & $-1$ & $-1$  & $Q_g$\;or\;$D$ & 32.036\\
   \rsp $(b_2)$ & $+1$ & $-1$ & $-1$ & $+1$ & $+1$ & $Q_g$\;or\;$D$ & 64.144\\
   \rsp $(b_3)$ & $-1$ & $+1$ & $-1$ & $+1$ & $+1$  & $Q_g$\;or\;$D$ & 64.144\\
   \rsp $(b_4)$ & $-1$ & $+1$ & $-1$ & $+1$ & $-1$ & $Q_g$\;or\;$D$ & 32.033\\
   \rsp $(b_5)$ & $+1$ & $-1$ & $+1$ & $-1$ & $-1$  & $D$ & 64.146\\
   \rsp $(b_6)$ & $+1$ & $-1$ & $-1$ & $+1$ & $-1$ & $D$ & 32.033\\
   \rsp $(b_7)$ & $+1$ & $-1$ & $+1$ & $+1$ & $-1$ & $D$ &  64.144\\
   \rsp $(b_8)$ & $-1$ & $+1$ & $+1$ & $+1$ & $-1$ & $S$\;or\;$D$ & 64.144 \\
   \rsp $(b_9)$ & $-1$ & $-1$ & $-1$ & $-1$ & $-1$ & $Q$ & 32.037 \\
   \rsp $(b_{10})$ & $-1$ & $-1$ & $-1$ & $+1$ & $+1$ & $Q$ & 64.147 \\
   \rsp $(b_{11})$ & $-1$ & $-1$ & $+1$ & $-1$ & $-1$  & $V_4$ & 32.036\\
   \rsp $(b_{12})$ & $-1$ & $-1$ & $-1$ & $+1$ & $-1$ & $V_4$ & 32.039 \\
   \rsp $(b_{13})$ & $-1$ & $-1$ & $+1$ & $+1$ & $-1$  & $V_4$ & 32.034\\
  \end{tabular}

\newpage

\noindent{\bf Type III}\;\;  ($d_1,d_2,d_3,d_4<0$, \;$d_j\not=-4$ ), \;\; $G\simeq V_4$

\medskip
$(d_1/p_2)=(d_1/p_4)=(d_4/p_3)=-1$ and

\medskip

\begin{tabular}{c|ccc|c}
   \text{case} & $(d_1/p_3)$ & $(d_2/p_3)$ & $(d_2/p_4)$ & $G^+/G_3^+$ \\ \hline
   \rsp $(c_1)$ & $+1$ & $-1$ & $+1$ & 32.036 \\
   \rsp $(c_2)$ & $+1$ & $+1$ & $+1$ & 32.033  \\
   \rsp $(c_3)$ & $-1$ & $+1$ & $-1$ & 64.150  \\
\end{tabular}

 \bigskip

 \noindent{\bf Type IV}\;\;  ($d_1,d_2,d_3,d_4<0$, \;$d_4=-4$ ),\;\; $G\simeq V_4$

\medskip

$(d_1/p_2)=(d_2/p_3)=-1$

\medskip
\begin{tabular}{c|cccc|c}
   \text{case} & $(d_1/p_3)$ & $(d_1/2)$ & $(d_2/2)$
            & $(d_3/2)$ & $G^+/G_3^+$ \\ \hline
   \rsp $(d_1)$ & $+1$ & $+1$ & $+1$ & $+1$ & 64.150\\
   \rsp $(d_2)$ & $+1$ & $+1$ & $-1$ & $+1$ & 32.036 \\
   \rsp $(d_3)$ & $+1$ & $-1$ & $-1$ & $+1$ & 32.037  \\
   \rsp $(d_4)$ & $+1$ & $-1$ & $-1$ & $-1$  & 32.041 \\
   \rsp $(d_5)$ & $-1$ & $+1$ & $-1$ & $+1$ & 32.033 \\
   \rsp $(d_6)$ & $-1$ & $+1$ & $+1$ & $-1$ & 32.036  \\
   \rsp $(d_7)$ & $-1$ & $-1$ & $+1$ & $-1$ & 32.036 \\
   \rsp $(d_8)$ & $-1$ & $-1$ & $-1$ & $+1$  & 32.033 \\
  \end{tabular}
\bigskip

The discriminant of $k$ is given by $d_k=d_1d_2d_3d_4$, where $d_i$ are prime discriminants as described  in the four types above. The cases give a complete list of those fields $k$ for which $\Cl_2(k)\simeq V_4$ and $d_k$ is not the sum of two squares. Furthermore, $G=\Gal(k^2/k)$, where $k^2$ is the second Hilbert $2$-class field of $k$, and $G^+=\Gal(k_+^2/k)$, where $k_+^2$ is the second narrow Hilbert $2$-class field of $k$.

\medskip
We have the following abbreviations:

\medskip

$V_4=(2,2)$, i.e., the Klein 4-group;

$Q$, the quaternion group of order $8$;

$Q_g$, the generalized quaternion group (of order $\geq 16$):

$D$, the dihedral group of order $2^m\geq 8$;

$S$, the semi-dihedral group (of order $m\geq 16$).
\medskip

Finally, the numbers associated with $G^+/G_3^+$ represent the $2$-groups  given in \cite{HS}.

\newpage

\section{{\bf Appendix~II}}

We have reproduced Appendix~II in \cite{BS3} for the convenience of the reader. A description of the tables is given at the end of the appendix.
\bigskip

{\tiny $$\begin{array}{|c|c|c|c|c|c|}\hline
\rsp {\rm Case} &  (\nu_{12},\nu_{13},\nu_{23},\nu_{14},\nu_{24},\nu_{34})&\delta\qquad \delta_1\qquad \delta_2 & N\varepsilon_{12}& \Cl_2(k) &\\
  & N_1\;\;\;\;N_2\;\;\;\;N_3 & \ker j_1\;\;\;\ker j_2\;\;\;\ker j_3  &  & \text{Type of $G$} & G^+/G^+_3  \\
  & q(k_1)\;\;q(k_2)\;\; q(k_3) & h_2(k_1)\;\;\;h_2(k_2)\;\;\;h_2(k_3)&  &  \text{Order of $G$} &\\\hline\hline
a_1 &(0,0,0,1,1,0) & p_1p_2p_4\;\;\;\;\;\;p_2p_4\;\;\;\;\;\;p_1p_4 & \left(
                                                             \begin{array}{c}
                                                               -1 \\
                                                               +1 \\
                                                             \end{array}
                                                           \right)      & \la [\fp_1],[\fp_2]\ra & \\
& \la [\fp_2]\ra\;\; \la [\fp_1]\ra \;\;\la [\fp_1\fp_2]\ra & \la [\fp_1]\ra \;\;\; \la [\fp_2]\ra \;\;\; \left(
                                                                                                 \begin{array}{c}
                                                                                                   \la [\fp_1\fp_2]\ra \\
                                                                                                   \la [\fp_1],[\fp_2]\ra   \\
                                                                                                 \end{array}
                                                                                               \right) & & \left(
                                                                                                             \begin{array}{c}
                                                                                                               Q_g \\
                                                                                                               D \\
                                                                                                             \end{array}
                                                                                                           \right)
                                                                                        & 64.144        \\
 & 2\qquad 2\qquad 2 & 4\qquad 4\qquad 2h_2(d_1d_2) & & 4h_2(d_1d_2) &\\\hline
a_2 &(0,1,1,0,0,0) & p_4\;\;\;\;\;\;p_4\;\;\;\;\;\;p_4 & \left(
                                                             \begin{array}{c}
                                                               -1 \\
                                                               +1 \\
                                                             \end{array}
                                                           \right)      & \la [\fp_1],[\fp_2]\ra  & \\
& \la [\fp_2]\ra\;\; \la [\fp_1]\ra \;\;\la [\fp_1\fp_2]\ra & \la [\fp_1]\ra \;\;\; \la [\fp_2]\ra \;\;\; \left(
                                                                                                 \begin{array}{c}
                                                                                                   \la [\fp_1\fp_2]\ra \\
                                                                                                   \la [\fp_1],[\fp_2]\ra   \\
                                                                                                 \end{array}
                                                                                               \right) & & \left(
                                                                                                             \begin{array}{c}
                                                                                                               Q_g \\
                                                                                                               D \\
                                                                                                             \end{array}
                                                                                                           \right)
                                                                                           & 64.144     \\
 & 2\qquad 2\qquad 2 & 4\qquad 4\qquad 2h_2(d_1d_2) & & 4h_2(d_1d_2) &\\\hline
a_3 &(1,0,1,0,0,0) & p_4 \;\;\; p_4 \:\:\: \left(
                                                             \begin{array}{c}
                                                               p_4 \\
                                                               p_1 \;\text{or}\;p_1p_4 \\
                                                             \end{array}
                                                           \right)     & -1   & \la [\fp_1],[\fp_2]\ra & \\
& \la [\fp_1\fp_2]\ra\;\; \la [\fp_2]\ra \;\;\la [\fp_1]\ra & \la [\fp_1]\ra \;\;\; \left(
                                                                \begin{array}{c}
                                                                  \la [\fp_2]\ra \\
                                                                   \la [\fp_1],[\fp_2]\ra \\
                                                                \end{array}
                                                              \right)
 \;\;\; \la [\fp_1\fp_2]\ra  & & \left(
                                                                                                             \begin{array}{c}
                                                                                                               Q_g \\
                                                                                                               D \\
                                                                                                             \end{array}
                                                                                                           \right)
                                                                                       & 64.144         \\
 & 2\quad \left(
     \begin{array}{c}
       2 \\
       1 \\
     \end{array}
   \right)
   \quad 2 & 4\quad \left(
                       \begin{array}{c}
                         2h_2(d_1d_3d_4) \\
                         h_2(d_1d_3d_4) \\
                       \end{array}
                     \right)
   \quad 4  & & \left(
                       \begin{array}{c}
                         4h_2(d_1d_3d_4) \\
                         2h_2(d_1d_3d_4) \\
                       \end{array}
                     \right) & \\\hline
 a_4 &(1,0,0,0,1,0) & p_1p_2p_4\;\;\; p_2p_4 \;\;\;\left(
                                                             \begin{array}{c}
                                                               p_1p_4 \\
                                                               p_1 \;\text{or}\;p_4 \\
                                                             \end{array}
                                                           \right)     & -1   & \la [\fp_1],[\fp_2]\ra & \\
& \la [\fp_1\fp_2]\ra\;\; \la [\fp_2]\ra \;\;\la [\fp_1]\ra & [\fp_1] \;\;\;\left(
                                                                \begin{array}{c}
                                                                  \la [\fp_2]\ra \\
                                                                   \la [\fp_1],[\fp_2]\ra \\
                                                                \end{array}
                                                              \right)
 \;\;\; \la [\fp_1\fp_2]\ra   & & \left(
                                                                                                             \begin{array}{c}
                                                                                                               Q_g \\
                                                                                                               D \\
                                                                                                             \end{array}
                                                                                                           \right)
                                                                                           &64.144     \\
 & 2 \quad \left(
     \begin{array}{c}
       2 \\
       1 \\
     \end{array}
   \right)
   \quad 2 & 4\quad  \left(
                       \begin{array}{c}
                         2h_2(d_1d_3d_4) \\
                         h_2(d_1d_3d_4) \\
                       \end{array}
                     \right)
   \quad 4 & & \left(
                       \begin{array}{c}
                         4h_2(d_1d_3d_4) \\
                         2h_2(d_1d_3d_4) \\
                       \end{array}
                     \right) & \\\hline
a_5 &(0,1,0,1,1,0) & p_2p_3p_4\;\;\;\;\;\;p_2p_4\;\;\;\;\;\;p_3p_4 & \left(
                                                             \begin{array}{c}
                                                               -1 \\
                                                               +1 \\
                                                             \end{array}
                                                           \right)      & \la [\fp_2],[\fp_3]\ra & \\
& \la [\fp_2]\ra\;\; \la [\fp_3]\ra \;\;\la [\fp_2\fp_3]\ra & \la [\fp_3]\ra \;\;\; \la [\fp_2]\ra \;\;\;                                                                                                 \la [\fp_2],[\fp_3]\ra    & & D &32.033 \\
 & 2\qquad 2\qquad \left(
                     \begin{array}{c}
                       1\\
                       2 \\
                     \end{array}
                   \right)
  & 4\qquad 4\qquad \left(
                                       \begin{array}{c}
                                         h_2(d_1d_2) \\
                                         2h_2(d_1d_2) \\
                                       \end{array}
                                     \right)
&  & \left(
                                       \begin{array}{c}
                                         2h_2(d_1d_2) \\
                                         4h_2(d_1d_2) \\
                                       \end{array}
                                     \right) &  \\\hline
a_6 &(0,1,1,0,1,0) & p_1p_3p_4\;\;\;\;\;\;p_3p_4\;\;\;\;\;\;p_4 & \left(
                                                             \begin{array}{c}
                                                               -1 \\
                                                               +1 \\
                                                             \end{array}
                                                           \right)      & \la [\fp_1],[\fp_3]\ra & \\
& \la [\fp_1\fp_3]\ra\;\; \la [\fp_1]\ra \;\;\la [\fp_3]\ra & \la [\fp_1]\ra \;\;\; \la [\fp_1\fp_3]\ra \;\;\;                                                                                                 \la [\fp_1],[\fp_3]\ra    & & D & 32.033 \\
 & 2\qquad 2\qquad \left(
                     \begin{array}{c}
                       1\\
                       2 \\
                     \end{array}
                   \right)
  & 4\qquad 4\qquad \left(
                                       \begin{array}{c}
                                         h_2(d_1d_2) \\
                                         2h_2(d_1d_2) \\
                                       \end{array}
                                     \right)
&  & \left(
                                       \begin{array}{c}
                                         2h_2(d_1d_2) \\
                                         4h_2(d_1d_2) \\
                                       \end{array}
                                     \right) & \\\hline
a_7 &(0,1,0,0,1,0) & p_2p_4\;\;\;\;\;\;p_2p_4\;\;\;\;\;\; p_4 & \left(
                                                             \begin{array}{c}
                                                               -1 \\
                                                               +1 \\
                                                             \end{array}
                                                           \right)      & \la [\fp_1],[\fp_2]\ra  &\\
& \la [\fp_2]\ra\;\; \la [\fp_1]\ra \;\;\la [\fp_1\fp_2]\ra & \la [\fp_1]\ra \;\;\; \la [\fp_2]\ra \;\;\;                                                                                                 \la [\fp_1],[\fp_2]\ra    & & D & 64.144\\
 & 2\qquad 2\qquad \left(
                     \begin{array}{c}
                       1\\
                       2 \\
                     \end{array}
                   \right)
  & 4\qquad 4\qquad \left(
                                       \begin{array}{c}
                                         h_2(d_1d_2) \\
                                         2h_2(d_1d_2) \\
                                       \end{array}
                                     \right)
&  & \left(
                                       \begin{array}{c}
                                         2h_2(d_1d_2) \\
                                         4h_2(d_1d_2) \\
                                       \end{array}
                                     \right) &  \\\hline

a_8 &(1,0,1,0,1,0) & p_3p_4\;\;\; p_3p_4 \;\;\; \left(
                                                             \begin{array}{c}
                                                               p_1 \\
                                                               p_4 \;\text{or}\;p_1p_4 \\
                                                             \end{array}
                                                           \right)    & -1   & \la [\fp_1],[\fp_3]\ra &\\
& \la [\fp_3]\ra\;\; \la [\fp_1\fp_3]\ra \;\;\la [\fp_1]\ra & \la [\fp_1]\ra \;\;\; \left(
                                                                \begin{array}{c}
                                                                  \la [\fp_1]\ra \\
                                                                   \la [\fp_1],[\fp_3]\ra \\
                                                                \end{array}
                                                              \right)
 \;\;\; \la [\fp_3]\ra  & & \left(
                                                                                                             \begin{array}{c}
                                                                                                               S \\
                                                                                                               D \\
                                                                                                             \end{array}
                                                                                                           \right)
                                                                                   & 32.033             \\
 & 2\quad \left(
     \begin{array}{c}
       2 \\
       1 \\
     \end{array}
   \right)
   \quad 2 & 4\quad \left(
                       \begin{array}{c}
                         2h_2(d_1d_3d_4) \\
                         h_2(d_1d_3d_4) \\
                       \end{array}
                     \right)
   \quad 4  & & \left(
                       \begin{array}{c}
                         4h_2(d_1d_3d_4) \\
                         2h_2(d_1d_3d_4) \\
                       \end{array}
                     \right) & \\\hline
\rsp a_9 &(1,0,0,1,1,0) & p_1p_2p_4\;\;\;\;\;\;p_2p_4\;\;\;\;\;\;p_1p_4 & -1     & \la [\fp_1],[\fp_2]\ra & \\
& \la [\fp_1]\ra\;\; \la [\fp_2]\ra \;\;\la [\fp_1\fp_2]\ra & \la [\fp_1]\ra \;\;\; \la [\fp_2]\ra \;\;\; \la [\fp_1\fp_2]\ra  & & Q & 64.147 \\
 & 2\qquad 2\qquad 2 & 4\qquad 4\qquad 4 & & 8 & \\\hline
\rsp a_{10} &(1,1,1,0,0,0) & p_4\;\;\;\;\;\;p_4\;\;\;\;\;\; p_4 & -1     & \la [\fp_1],[\fp_2]\ra &  \\
& \la [\fp_1]\ra\;\; \la [\fp_2]\ra \;\;\la [\fp_1\fp_2]\ra & \la [\fp_1]\ra \;\;\; \la [\fp_2]\ra \;\;\; \la [\fp_1\fp_2]\ra  & & Q & 64.147 \\
 & 2\qquad 2\qquad 2 & 4\qquad 4\qquad 4 & & 8 &\\\hline
\rsp a_{11} &(1,1,0,1,1,0) & p_2p_3\;\;\;\;\;\;p_2p_4\;\;\;\;\;\; p_3p_4 & -1     & \la [\fp_1],[\fp_2]\ra & \\
& \la [\fp_1\fp_2]\ra\;\; \la [\fp_2]\ra \;\;\la [\fp_1]\ra & \la [\fp_1],[\fp_2]\ra \;\;\; \la [\fp_1],[\fp_2]\ra \;\;\; \la [\fp_1],[\fp_2]\ra  & & (2,2) & 32.39  \\
 & 1\qquad 1\qquad 1 & 2\qquad 2\qquad 2 & & 4 & \\\hline
\rsp a_{12} &(1,1,1,0,1,0) & p_2p_3\;\;\;\;\;\;p_3p_4\;\;\;\;\;\; p_4 & -1     & \la [\fp_1],[\fp_2]\ra & \\
& \la [\fp_1]\ra\;\; \la [\fp_1\fp_2]\ra \;\;\la [\fp_2]\ra & \la [\fp_1],[\fp_2]\ra \;\;\; \la [\fp_1],[\fp_2]\ra \;\;\; \la [\fp_1],[\fp_2]\ra  & & (2,2) & 32.39  \\
 & 1\qquad 1\qquad 1 & 2\qquad 2\qquad 2 & & 4 &\\\hline
\rsp a_{13} &(1,1,0,0,1,0) & p_2p_3\;\;\;\;\;\;p_2p_4\;\;\;\;\;\; p_4 & -1     & \la [\fp_1],[\fp_2]\ra & \\
& \la [\fp_1]\ra\;\; \la [\fp_2]\ra \;\;\la [\fp_1\fp_2]\ra & \la [\fp_1],[\fp_2]\ra \;\;\; \la [\fp_1],[\fp_2]\ra \;\;\; \la [\fp_1],[\fp_2]\ra  & & (2,2) & 32.34 \\
 & 1\qquad 1\qquad 1 & 2\qquad 2\qquad 2 & & 4 &\\\hline
\end{array}$$  }

{\tiny  $$\begin{array}{|c|c|c|c|c|c|}\hline
\rsp {\rm Case} & (\nu_{12},\nu_{13},\nu_{23},\nu_{14},\nu_{24},\nu_{34})&\delta\qquad \delta_1\qquad \delta_2 & N\varepsilon_{12}& \Cl_2(k) &\\
  & N_1\;\;\;\;N_2\;\;\;\;N_3 & \ker j_1\;\;\;\ker j_2\;\;\;\ker j_3  &  & \text{Type of $G$} & G^+/G^+_3  \\
  & q(k_1)\;\;q(k_2)\;\; q(k_3) & h_2(k_1)\;\;\;h_2(k_2)\;\;\;h_2(k_3)&  &  \text{Order of $G$} &\\\hline\hline
 b_1 &(0,1,1,1,1,\left[
                              \begin{array}{c}
                                1 \\
                                0 \\
                              \end{array}
                            \right]
 ) & \left[
                              \begin{array}{c}
                                2p_3 \\
                               2p_1p_2 \\
                              \end{array}
                            \right]\;\;\;\;\;\; \left[
                              \begin{array}{c}
                                2p_2 \\
                               2p_3 \\
                              \end{array}
                            \right]\;\;\;\;\;\;  \left[
                              \begin{array}{c}
                                2p_1 \\
                               2p_3 \\
                              \end{array}
                            \right] & \left(
                                                             \begin{array}{c}
                                                               -1 \\
                                                               +1 \\
                                                             \end{array}
                                                           \right)      & \la [\fp_1],[\fp_2]\ra & \\
& \la [\fp_2]\ra\;\; \la [\fp_1]\ra \;\;\la [\fp_1\fp_2]\ra & \la [\fp_1]\ra \;\;\; \la [\fp_2]\ra \;\;\left(
                                                                                                         \begin{array}{c}
                                                                                                           \la [\fp_1\fp_2]\ra \\
                                                                                                           \la [\fp_1],[\fp_2]\ra \\
                                                                                                         \end{array}
                                                                                                       \right)
 & & \left(
                                                             \begin{array}{c}
                                                               Q_g \\
                                                               D \\
                                                             \end{array}
                                                           \right)
                                                                                        & 32.36        \\
 & 2\qquad 2\qquad 2  & 4\qquad 4\qquad 2h_2(p_1p_2)  & & 4h_2(p_1p_2)  & \\\hline
b_2 &(0,1,1,0,0,\left[
                              \begin{array}{c}
                                1 \\
                                0 \\
                              \end{array}
                            \right]) & \left[
                              \begin{array}{c}
                                2p_1p_2p_3 \\
                                     2\\
                              \end{array}
                            \right]\;\;\;\;\;\; \left[
                              \begin{array}{c}
                                2p_2p_3 \\
                               2 \\
                              \end{array}
                            \right]\;\;\;\;\;\;  \left[
                              \begin{array}{c}
                                2p_1p_3 \\
                               2 \\
                              \end{array}
                            \right] & \left(
                                                             \begin{array}{c}
                                                               -1 \\
                                                               +1 \\
                                                             \end{array}
                                                           \right)      & \la [\fp_1],[\fp_2]\ra & \\
& \la [\fp_2]\ra\;\; \la [\fp_1]\ra \;\;\la [\fp_1\fp_2]\ra & \la [\fp_1]\ra \;\;\; \la [\fp_2]\ra \;\;\left(
                                                                                                         \begin{array}{c}
                                                                                                           \la [\fp_1\fp_2]\ra \\
                                                                                                           \la [\fp_1],[\fp_2]\ra \\
                                                                                                         \end{array}
                                                                                                       \right)
 & & \left(
                                                             \begin{array}{c}
                                                               Q_g \\
                                                               D \\
                                                             \end{array}
                                                           \right)
                                                                                        & 64.144        \\
 & 2\qquad 2\qquad 2  & 4\qquad 4\qquad 2h_2(p_1p_2)  & & 4h_2(p_1p_2)  & \\\hline
b_{3} &(1,0,1,0,0,\left[
                              \begin{array}{c}
                                1 \\
                                0 \\
                              \end{array}
                            \right]) & \left[
                              \begin{array}{c}
                                2p_1p_2p_3 \\
                                     2\\
                              \end{array}
                            \right]\;\;\;\;\;\; \left[
                              \begin{array}{c}
                                2p_2p_3 \\
                               2 \\
                              \end{array}
                            \right]\;\;\;\;\;\;   \;\;\; \left[\begin{array}{c}
                            \left(
                                                 \begin{array}{c}
                                                   2p_1p_3 \\
                                                   p_1\;\text{or}\;2p_3 \\
                                                 \end{array}
                                               \right)\\
                                                \left(
                                                 \begin{array}{c}
                                                   2 \\
                                                   p_1\;\text{or}\;2p_1 \\
                                                 \end{array}
                                               \right)\end{array}\right]
  & -1    & \la [\fp_1],[\fp_2]\ra & \\
& \la [\fp_1\fp_2]\ra\;\; \la [\fp_2]\ra \;\;\la [\fp_1]\ra & \la [\fp_1]\ra \;\;\left(
                                                                \begin{array}{c}
                                                                  \la [\fp_2]\ra \\
                                                                  \la [\fp_1],[\fp_2]\ra \\
                                                                \end{array}
                                                              \right)
 \; \;\; \la [\fp_1\fp_2]\ra
 & & \left(
        \begin{array}{c}
          Q_g \\
          D \\
        \end{array}
      \right)
                                                                                        & 64.144     \\
 & 2\qquad \left(
             \begin{array}{c}
               2 \\
               1 \\
             \end{array}
           \right)
 \qquad 2  & 4\qquad \left(
                       \begin{array}{c}
                         2h_2(p_1p_3) \\
                         h_2(p_1p_3) \\
                       \end{array}
                     \right)
 \qquad 4  & & \left(
                 \begin{array}{c}
                   4h_2(p_1p_3) \\
                   2h_2(p_1p_3) \\
                 \end{array}
               \right)
   & \\\hline
b_4 &(1,0,1,0,1,\left[
                              \begin{array}{c}
                                1 \\
                                0 \\
                              \end{array}
                            \right]) & \left[
                              \begin{array}{c}
                                2p_1p_2 \\
                               2p_3 \\
                              \end{array}
                            \right]\;\;\;\;\;\; \left[
                              \begin{array}{c}
                                2p_2 \\
                               2p_3 \\
                              \end{array}
                            \right] \;\;\; \left[\begin{array}{c}
                            \left(
                                                 \begin{array}{c}
                                                   2p_3 \\
                                                  p_1\;\text{or}\;2p_1p_3 \\
                                                 \end{array}
                                               \right)\\
                                                \left(
                                                 \begin{array}{c}
                                                   2p_1 \\
                                                   2\;\text{or}\;p_1 \\
                                                 \end{array}
                                               \right)\end{array}\right] & -1    & \la [\fp_1],[\fp_2]\ra & \\
& \la [\fp_1\fp_2]\ra\;\; \la [\fp_2]\ra \;\;\la [\fp_1]\ra & \la [\fp_1]\ra \; \;\;\left(
                                                                                                         \begin{array}{c}
                                                                                                           \la [\fp_2]\ra \\
                                                                                                           \la [\fp_1],[\fp_2]\ra \\
                                                                                                         \end{array}
                                                                                                       \right)\;\; \la [\fp_1\fp_2]\ra
 & & \left(
                                                             \begin{array}{c}
                                                               Q_g \\
                                                               D \\
                                                             \end{array}
                                                           \right)
                                                                                        & 32.033     \\
 & 2\qquad \left(
             \begin{array}{c}
               2 \\
               1 \\
             \end{array}
           \right)
 \qquad 2  & 4\qquad \left(
                       \begin{array}{c}
                         2h_2(p_1p_3) \\
                         h_2(p_1p_3) \\
                       \end{array}
                     \right)
 \qquad 4  & & \left(
                 \begin{array}{c}
                   4h_2(p_1p_3) \\
                   2h_2(p_1p_3) \\
                 \end{array}
               \right)
   & \\\hline
 b_5 &(0,1,0,1,1,\left[
                              \begin{array}{c}
                                1 \\
                                0 \\
                              \end{array}
                            \right]) & p_2\;\;\;\;\;\; p_2\;\;\;\;\;\; \left[
                              \begin{array}{c}
                                2p_1 \\
                                     2p_3  \\
                              \end{array}
                            \right] & \left(
                                                             \begin{array}{c}
                                                               -1 \\
                                                               +1 \\
                                                             \end{array}
                                                           \right)      & \la [\fp_1],[\fq]\ra & \\
& \la [\fp_1\fq]\ra\;\; \la [\fp_1]\ra \;\;\la [\fq]\ra & \la [\fp_1]\ra \;\;\; \la [\fp_1\fq]\ra \;\;\;
                                                                                                   \la [\fp_1],[\fq]\ra  & & D
                                                                                        & 64.146        \\
 & 2\qquad 2\qquad \left(\begin{array}{c}
                                                               1 \\
                                                               2 \\
                                                             \end{array}
                                                           \right) & 4\qquad 4\qquad \left(
                                                             \begin{array}{c}
                                                               h_2(p_1p_2) \\
                                                               2h_2(p_1p_2) \\
                                                             \end{array}
                                                           \right)  & & \left( \begin{array}{c}
                                                               2h_2(p_1p_2) \\
                                                               4h_2(p_1p_2) \\
                                                             \end{array} \right) & \\\hline

b_6 &(0,1,1,0,1,\left[
                              \begin{array}{c}
                                1 \\
                                0 \\
                              \end{array}
                            \right]) &  \left[
                              \begin{array}{c}
                                2p_2 \\
                               2p_1p_3 \\
                              \end{array}
                            \right]\;\;\;\;\;\; \left[
                              \begin{array}{c}
                                2p_2 \\
                               2p_3 \\
                              \end{array}
                            \right]\;\;\;\;\;\;  \left[
                              \begin{array}{c}
                                2p_1p_3 \\
                               2 \\
                              \end{array}
                            \right] & \left(
                                                             \begin{array}{c}
                                                               -1 \\
                                                               +1 \\
                                                             \end{array}
                                                           \right)      & \la [\fp_1],[\fp_2]\ra & \\
& \la [\fp_2]\ra\;\; \la [\fp_1]\ra \;\;\la [\fp_1\fp_2]\ra & \la [\fp_1]\ra \;\;\; \la [\fp_2]\ra \;\;\;
                                                                                                   \la [\fp_1],[\fp_2]\ra  & & D
                                                                                        & 32.033    \\
& 2\qquad 2\qquad \left(\begin{array}{c}
                                                               1 \\
                                                               2 \\
                                                             \end{array}
                                                           \right) & 4\qquad 4\qquad \left(
                                                             \begin{array}{c}
                                                               h_2(p_1p_2) \\
                                                               2h_2(p_1p_2) \\
                                                             \end{array}
                                                           \right)  & & \left( \begin{array}{c}
                                                               2h_2(p_1p_2) \\
                                                               4h_2(p_1p_2) \\
                                                             \end{array}\right)  & \\\hline

b_7 &(0,1,0,0,1,\left[
                              \begin{array}{c}
                                1 \\
                                0 \\
                              \end{array}
                            \right]) & p_2\;\;\;\;\;\; p_2\;\;\;\;\;\; \left[
                              \begin{array}{c}
                                2p_1p_3 \\
                                     2  \\
                              \end{array}
                            \right] & \left(
                                                             \begin{array}{c}
                                                               -1 \\
                                                               +1 \\
                                                             \end{array}
                                                           \right)      & \la [\fp_1],[\fq]\ra & \\
& \la [\fq]\ra\;\; \la [\fp_1]\ra \;\;\la [\fq\fp_1]\ra & \la [\fp_1]\ra \;\;\; \la [\fq]\ra \;\;\;
                                                                                                   \la [\fq],[\fp_1]\ra  & & D
                                                                                        & 64.144        \\
 & 2\qquad 2\qquad \left(\begin{array}{c}
                                                               1 \\
                                                               2 \\
                                                             \end{array}
                                                           \right) & 4\qquad 4\qquad \left(
                                                             \begin{array}{c}
                                                               h_2(p_1p_2) \\
                                                               2h_2(p_1p_2) \\
                                                             \end{array}
                                                           \right)  & & \left( \begin{array}{c}
                                                               2h_2(p_1p_2) \\
                                                               4h_2(p_1p_2) \\
                                                             \end{array}\right)  & \\\hline

b_8 &(1,0,0,0,1,\left[
                              \begin{array}{c}
                                1 \\
                                0 \\
                              \end{array}
                            \right]) & p_1p_2\;\;\;\;\;\;  p_2 \;\;\; \left[\begin{array}{c}
                            \left(
                                                 \begin{array}{c}
                                                   p_1 \\
                                                  2p_3\;\text{or}\;2p_1p_3 \\
                                                 \end{array}
                                               \right)\\
                                                \left(
                                                 \begin{array}{c}
                                                   p_1 \\
                                                   2\;\text{or}\;2p_1 \\
                                                 \end{array}
                                               \right)\end{array}\right]  & -1    & \la [\fp_1],[\fq]\ra & \\
& \la [\fq]\ra\;\; \la [\fq\fp_1]\ra \;\;\la [\fp_1]\ra & \la [\fp_1]\ra \; \;\;\left(
                                                                                                         \begin{array}{c}
                                                                                                           \la [\fp_1]\ra \\
                                                                                                           \la [\fq],[\fp_1]\ra \\
                                                                                                         \end{array}
                                                                                                       \right)\;\; \la [\fq]\ra
 & & \left(
                                                             \begin{array}{c}
                                                               S \\
                                                               D \\
                                                             \end{array}
                                                           \right)
                                                                                        & 64.144        \\
 & 2\qquad \left(
             \begin{array}{c}
               2 \\
               1 \\
             \end{array}
           \right)
 \qquad 2  & 4\qquad \left(
                       \begin{array}{c}
                         2h_2(p_1p_3) \\
                         h_2(p_1p_3) \\
                       \end{array}
                     \right)
 \qquad 4  & & \left(
                 \begin{array}{c}
                   4h_2(p_1p_3) \\
                   2h_2(p_1p_3) \\
                 \end{array}
               \right)
   & \\\hline

\rsp  b_{9} &(1,1,1,1,1,\left[
                              \begin{array}{c}
                                1 \\
                                0 \\
                              \end{array}
                            \right]) &  \left[
                              \begin{array}{c}
                                2p_3 \\
                               2p_1p_2 \\
                              \end{array}
                            \right]\;\;\;\;\;\; \left[
                              \begin{array}{c}
                                2p_2 \\
                               2p_3 \\
                              \end{array}
                            \right]\;\;\;\;\;\;  \left[
                              \begin{array}{c}
                                2p_1 \\
                               2p_3 \\
                              \end{array}
                            \right]  & -1    & \la [\fp_1],[\fp_2]\ra & \\
& \la [\fp_1]\ra\;\; \la [\fp_2]\ra \;\;\la [\fp_1\fp_2]\ra & \la [\fp_1]\ra \; \;\;\la [\fp_2]\ra \;\; \la [\fp_1\fp_2]\ra
 & & Q
                                                                                        & 32.037     \\
 & 2\qquad 2
 \qquad 2  & 4\qquad 4
 \qquad 4  & & 8
   & \\\hline
\rsp b_{10} &(1,1,1,0,0,\left[
                              \begin{array}{c}
                                1 \\
                                0 \\
                              \end{array}
                            \right]) & \left[
                              \begin{array}{c}
                                2p_1p_2p_3 \\
                               2 \\
                              \end{array}
                            \right]\;\;\;\;\;\; \left[
                              \begin{array}{c}
                                2p_2p_3 \\
                               2 \\
                              \end{array}
                            \right]\;\;\;\;\;\;  \left[
                              \begin{array}{c}
                                2p_1p_3 \\
                               2 \\
                              \end{array}
                            \right]  & -1    & \la [\fp_1],[\fp_2]\ra & \\
& \la [\fp_1]\ra\;\; \la [\fp_2]\ra \;\;\la [\fp_1\fp_2]\ra & \la [\fp_1]\ra \; \;\;\la [\fp_2]\ra \;\; \la [\fp_1\fp_2]\ra
 & & Q
                                                                                        & 64.147     \\
 & 2\qquad 2
 \qquad 2  & 4\qquad 4
 \qquad 4  & & 8
   & \\\hline
\rsp b_{11} &(1,1,0,1,1,\left[
                              \begin{array}{c}
                                1 \\
                                0 \\
                              \end{array}
                            \right]) &  \left[
                              \begin{array}{c}
                                2p_1p_3 \\
                               2p_2 \\
                              \end{array}
                            \right]\;\;\;\;\;\; p_2\;\;\;\;\;\;  \left[
                              \begin{array}{c}
                                2p_1 \\
                               2p_3 \\
                              \end{array}
                            \right]  & -1    & \la [\fp_1],[\fp_2]\ra & \\
& \la [\fp_1]\ra\;\; \la [\fp_1\fp_2]\ra \;\;\la [\fp_2]\ra & \la [\fp_1],[\fp_2]\ra \; \;\;\la [\fp_1],[\fp_2]\ra \;\; \la [\fp_1],[\fp_2]\ra
 & & (2,2)
                                                                                        & 32.036     \\
 & 1\qquad 1
 \qquad 1  & 2\qquad 2
 \qquad 2  & & 4
   & \\\hline
\rsp b_{12} &(1,1,1,0,1,\left[
                              \begin{array}{c}
                                1 \\
                                0 \\
                              \end{array}
                            \right]) & \left[
                              \begin{array}{c}
                                2p_1 \\
                               2p_2p_3 \\
                              \end{array}
                            \right]\;\;\;\;\;\; \left[
                              \begin{array}{c}
                                2p_2 \\
                               2p_3 \\
                              \end{array}
                            \right]\;\;\;\;\;\;  \left[
                              \begin{array}{c}
                                2p_1p_3 \\
                               2 \\
                              \end{array}
                            \right]  & -1    & \la [\fp_1],[\fp_2]\ra & \\
& \la [\fp_1]\ra\;\; \la [\fp_2]\ra \;\;\la [\fp_1\fp_2]\ra & \la [\fp_1],[\fp_2]\ra \; \;\;\la [\fp_1],[\fp_2]\ra \;\; \la [\fp_1],[\fp_2]\ra
 & & (2,2)
                                                                                        & 32.039     \\
 & 1\qquad 1
 \qquad 1  & 2\qquad 2
 \qquad 2  & & 4
   & \\\hline
\rsp b_{13} &(1,1,0,0,1,\left[
                              \begin{array}{c}
                                1 \\
                                0 \\
                              \end{array}
                            \right]) &  \left[
                              \begin{array}{c}
                                2p_1 \\
                               2p_2p_3 \\
                              \end{array}
                            \right]\;\;\;\;\;\; p_2\;\;\;\;\;\;  \left[
                              \begin{array}{c}
                                2p_1p_3 \\
                               2 \\
                              \end{array}
                            \right]   & -1    & \la [\fp_1],[\fp_2]\ra & \\
& \la [\fp_1]\ra\;\; \la [\fp_1\fp_2]\ra \;\;\la [\fp_2]\ra & \la [\fp_1],[\fp_2]\ra \; \;\;\la [\fp_1],[\fp_2]\ra \;\; \la [\fp_1],[\fp_2]\ra
 & & (2,2)
                                                                                        & 32.034     \\
 & 1\qquad 1
 \qquad 1  & 2\qquad 2
 \qquad 2  & & 4
   & \\\hline

\end{array}$$}

\newpage

 Here's a description of the tables above.

Column 1: This column indicates the cases $(a_i)$ and $(b_i)$ as given in Appendix~I.

For the moment, we skip to Column 4.

Column 4: Shows the possible values of $N\eps_{12}$, the norm of the fundamental unit $\eps_{12}$ in $\Q(\sqrt{d_1d_2}\,)$.
When there are two possible values in any entry, we represent them in a $2$-d column vector given in rounded brackets. If there are several such column vectors in a given case, all their rows correspond to the same relevant data.

Column 2: Row 1 gives the relevant values $\nu_{ij}\in\{0,1\}$ such that $(-1)^{\nu_{ij}}=(d_i/p_j)$, where $p_i$ is the rational prime dividing $d_i$. Moreover, in the cases $(b_i)$ we include the two possible values of $\nu_{34}$, which is given as a column vector in square brackets. (Observe here that the only changes occurring are the values of $\delta$, $\delta_1$, and $\delta_2$ in Column~3.) Row 2 gives $N_i=N_{k_i/k}\Cl_2(k_i)$ for $k_1=k(\sqrt{d_1}\,)$, $k_2=k(\sqrt{d_2}\,)$, and $k_3=k(\sqrt{d_1d_2}\,)$.  Row 3 gives the unit indices  $q(k_i)=q(k_i/\Q)$, which is the index of $E_{k_i}$ modulo the subgroup generated by the unit groups in the quadratic subfields of $k_i$. (Computing these indices requires  data in the first row of the third column.)

Column 3: Row 1 gives the values $\delta=\delta(\varepsilon_k)=\sfk N(1+\varepsilon_k)$, where $\sfk$ denotes the square-free kernel. Similarly, $\delta_1=\delta(\eps_{234})$ and $\delta_2=\delta(\eps_{134})$. Here $\eps_{i34}$ denotes the fundamental unit of $\Q(\sqrt{d_id_3d_4}\,)$. We also observe that $\delta_3=\delta(\eps_{34})=p_4$ in the $(a_i)$ cases and in the $(b_i)$ cases $\delta(\eps_{34})=2p_3$ or $2$, according as $p_3\equiv 3$ or $7\bmod 8$, respectively, cf.\,\cite{Kub}.  Rows 2 and 3 are self-explanatory.
(Here $h_2(k_i)$ denotes the $2$-class number of $k_i$.)

Column 5: Is self-explanatory.

Column 6: Gives the quotient group $G^+/G_3^+$ for  $G^+=\Gal(k_+^2/k)$.

Finally, for the $(b_i)$ cases, $\fq$ denotes the prime ideal in $k$ above $2$.

\bigskip

{\bf Addresses}

\medskip

\author{Elliot Benjamin}

\address{School of Social and Behavioral Sciences,  Capella University,

  Minneapolis, MN 55402, USA; www.capella.edu}

\email{ben496@prexar.com}

\medskip

\author{C.\,Snyder}

\address{Department of Mathematics and Statistics, UMaine, Orono, ME 04469, USA}

\email{wsnyder@maine.edu}

\end{document}